\documentclass[10pt,a4paper,oneside]{article}

\usepackage{amssymb}
\usepackage{amsmath}
\usepackage{hyperref}
\usepackage{stmaryrd}
\usepackage{slashed}
\usepackage{url}
\usepackage{rotating}

\newcounter{count}[section]
\renewcommand{\thecount}{\arabic{section}.\arabic{count}}

\numberwithin{equation}{section}

\newenvironment{Umgeb1rekursiv}[1]
   {\vspace{0.5cm}
     \noindent
     \par\noindent
     \refstepcounter{count}
     \textbf{#1~\thecount}
     \hspace{0.2cm}
     \itshape
   }
   {\vspace{0.2cm}\par}

\newenvironment{Umgeb1}[1]
   {\vspace{0.5cm}
     \par\noindent
     \refstepcounter{count}
     \textbf{#1~\thecount}
     \hspace{0.2cm}
   }
   {\vspace{0.2cm}\par}

\newenvironment{UmgebBeweis}[1]
   {\vspace{0.5cm}
     \par\noindent
     \textbf{Proof.}
     \hspace{0.2cm}
   }
   {\hfill $\square$
   \vspace{0.2cm}}
   
\newenvironment{bsp}{\begin{Umgeb1}{Example}}
   {\end{Umgeb1}}
\newenvironment{sat}{\begin{Umgeb1rekursiv}{Proposition}}
  {\end{Umgeb1rekursiv}}
\newenvironment{lem}{\begin{Umgeb1rekursiv}{Lemma}}
  {\end{Umgeb1rekursiv}}
\newenvironment{kor}{\begin{Umgeb1rekursiv}{Corollary}}
  {\end{Umgeb1rekursiv}}  
\newenvironment{theorem}{\begin{Umgeb1rekursiv}{Theorem}}
  {\end{Umgeb1rekursiv}}  

\newenvironment{defn}{\begin{Umgeb1rekursiv}{Definition}}
  {\end{Umgeb1rekursiv}}
\newenvironment{bew}{\begin{UmgebBeweis}{Proof}}
   {\end{UmgebBeweis}}

\newenvironment{bem}{\begin{Umgeb1}{Remark}}
   {\end{Umgeb1}}
   
\newcommand{\R}{\mathbb{R}}
\newcommand{\Ham}{\mathbb{H}}
\newcommand{\C}{\mathbb{C}}
\newcommand{\N}{\mathbb{N}}

\newcommand{\id}{id}

\newcommand{\diver}{\operatorname{div}}

\newcommand{\tr}{tr}

\newcommand{\abs}[1]{\lvert#1\rvert}
\renewcommand{\part}[2]{\frac{\partial #1}{\partial #2}}

\newcommand{\grad}{\operatorname{grad}}

\title{Spinor Residue Family Operators and Spectral Theory of Dirac Operator for Poincar\'e-Einstein Metrics}
\author{Matthias Fischmann, Petr Somberg}
\begin{document}
\date{}
\maketitle
\begin{abstract}
We study conformal $Spin$-subgeometry of submanifolds in a semi-Riemannian $Spin$-manifold, focusing
on conformal $Spin$-manifolds $(M,[h])$ and their Poincar\'e-Einstein metrics $(X,g_+)$. Our approach is
based on the spectral theory of Dirac operator in the ambient $Spin$-manifold, and associated spinor 
valued meromorphic family of distributions with residues given by the residue family 
operators $\slashed{D}_N^{res}(h;\lambda)$ on spinors. We develop basic aspects and properties of 
$\slashed{D}_N^{res}(h;\lambda)$ including conformal covariance, factorization properties 
by conformally covariant operators for both flat and curved semi-Riemannian $Spin$-manifolds, 
and Poisson transformation.      
\end{abstract}

\vspace{0.3cm}					

{\bf Keywords:} Conformal semi-Riemannian $Spin$-geometry and subgeometry, Spectral theory of Dirac operator, 
Invariant distributions, Poisson transform, Conformal powers of the Dirac operator.

\vspace{0.1cm}

{\bf MSC2010:} 58C40, 58J50, 53A30, 53A55, 46F12.                

\allowdisplaybreaks


\section{Introduction}

The fact that the orthogonal Lie group may be regarded as the isometry group 
of a space and at the same time as the conformal group of another space has 
far-reaching consequences in geometrical analysis, representation theory and 
topology of manifolds. 

The representative example of this phenomenon is given by conformally equivariant operators 
for principal series representations on the conformal sphere $S^n$, which can 
be interpreted 
as scattering operators on the hyperbolic space; the Poisson transform allows 
to realize tensor-spinors on compactified boundary (the conformal sphere) of the hyperbolic 
space as asymptotics of eigenspaces for the algebra of invariant differential 
operators in the interior of conformal compactification, cf. \cite{Helgason}.
  
A curved generalization of the previous example is called ambient metric, see \cite{FG3}. 
It associates to a conformal manifold $(M,[h])$ of 
dimension $n\geq 3$ a pseudo-Riemannian Ricci flat ambient manifold of dimension $n+2$, so that conformal 
invariants of $(M,[h])$ are induced by pseudo-Riemannian invariants of the ambient metric.  
In particular, the ambient metric allows to construct conformally covariant operators
including GJMS operators as conformal modifications of  
powers of the Laplace operator, \cite{GJMS, GZ, GoverPeterson}, and conformal powers of the
Dirac operator, \cite{HS,GMP1,Fischmann}.  
In the ambient space there is a Lorentzian hypersurface $X=M\times(0,\varepsilon)$, and the 
ambient metric induces on $X$ an Einstein metric $g_+$ termed 
{\sf Poincar\'e-Einstein metric}. Thus $(M,[h])$ 
is realized as the conformal infinity of the Poincar\'e-Einstein metric $(X,g_+)$ and  
the geometric scattering theory for the Poincar\'e-Einstein metric produces the scattering 
operator $S(h,\lambda)$, \cite{GZ,GMP}, 
fulfilling covariant transformation property 
in the conformal class $[h]$. 

Beside the class of conformally covariant operators acting on sheaves over the same base manifold $M$ 
like the scattering or GJMS operators mentioned above, there is a sequence of $1$-parameter families of 
differential operators depending on $\lambda\in{\mathbb C}$ and called 
residue families,
\begin{align*}
   D_N^{res}(h,\lambda): C^\infty\big(M\times [0,\varepsilon)\big)\to C^\infty(M),
\end{align*}
$N\in\N_0$ ($0\leq N\leq n$ for even $n$), for $n$-dimensional conformal manifold $(M,[h])$, cf. \cite{Juhl}. 
The residue families exhibit a covariant transformation property
for conformal change of the metric $h$, specialize to GJMS operators at specific values of $\lambda$ 
and obstruct the existence of a continuation of the function $r^\mu u$,
$u\in \ker\big(\Delta_{g_+}+\lambda(n-\lambda)\big)$, $\mu\in{\mathbb C}$, to a distribution on 
$M\times [0,\varepsilon)$. The residue families encode neat invariants of conformal structure 
$[h]$ on $M$, e.g., the Branson's $Q$-curvature is produced 
as a derivative in the variable $\lambda$ of the critical residue families $D_n^{res}(h;\lambda)$, for even $n$, 
at the value $\lambda=0$, evaluated on the constant function $1$.
Furthermore, the residue families enjoy a system of factorization identities, given by 
pre-compositions resp. post-compositions with GJMS operators on $(M,h)$ resp. on the conformal 
compactification  
$(M\times [0,\varepsilon),\bar{g})$ of $(X,g_+)$. In the flat case, 
the factorization identities are the consequences 
of factorization properties in the representation theory, describing compositions of homomorphisms of 
generalized Verma modules, see \cite{KOSS}. To summarize, the residue families are differential 
invariants of conformal submanifolds and can be analytically realized in residues of a distribution
constructed out of the defining function of the submanifold and an eigenfunction of the Laplace operator
in the ambient manifold. 

The main results of our article are the construction and basic properties of spinor analogues of residue families,
thereby producing differential invariants of conformal $Spin$-submanifolds. In this case are the 
basic building blocks of residue family operators on spinors the conformally covariant operators 
with leading terms given by odd powers of the Dirac operator ${\slashed D}$ on $(M,h)$, 
cf. \cite{HS, GMP1, Fischmann}. We focus again on the case of Poincar\'e-Einstein metric 
and its conformal compactification. 

This theme conceptually fits into the framework of boundary valued 
problems and Poisson transform for the systems of partial 
differential equations with regular singularities along a submanifold, cf. \cite{KKMOOT, Helgason}.
A closely parallel topics include analytic aspects of geometric scattering theory 
for spinors on Poincar\'e-Einstein metrics $(\overline{X},g_+)$, cf. \cite{GMP,GMP1}, or algebraic aspects 
of the residue family operators on spinors based on the classification of homomorphisms of generalized 
Verma modules for conformal parabolic subalgebras and inducing spinor representations, cf. \cite{KOSS}.   

\vspace{5pt}

Let us briefly review the content of the present article. We start with 
the formulation of boundary valued problem for spinors, focusing on the 
eigenvalue equation for the Dirac 
operator on Poincar\' e-Einstein metric associated to 
$(M,[h])$, cf. Section \ref{Boundary}. 
In Section \ref{ResidueFamilies}, 
we construct the residue family operators on spinors, cf. Definition \ref{DefResFamOp},
through the analysis of obstructions for extending formal asymptotic eigenspinors as spinor-valued distributions, 
supported on $M\times [0,\varepsilon)$ and acting on the functional space of compactly supported smooth spinors. 
We discuss the general case of curved semi-Riemannian $Spin$-manifold in a way that for the flat 
case the residue family operators correspond to conformally covariant differential operators induced by homomorphisms
of Verma modules for codimension one orthogonal Lie algebras and their conformal parabolic 
subalgebras, cf. Theorem \ref{Intertwiner}.
Section \ref{ConformalTrafo} contains a short digression on the conformal covariance of residue family operators, 
cf. Theorem \ref{ConformalTrafoLawResidueFamily}. In addition, we define a family of first order differential 
operators associated to arbitrary hypersurface in a conformal manifold and prove its conformal covariance, 
cf. Proposition \ref{GeneralHypersurfaceCovariance}.  
In the flat case, we determine the 
full system of factorization identities for residue family operators on spinors, cf. Theorems \ref{FactorizationSmallDirac} and 
\ref{FactorziationBicDirac}, in Section \ref{FactorizationProp}. 
The abstract conclusion is the consequence
of the classification of conformally covariant differential operators on spinors, and the proof is based on combinatorial 
identities for the hypergeometric functions. We also prove factorization identities for a few residue 
family operators in the general curved case, exploiting the metric construction of generalized cylinders over
conformal boundary $(M,[h])$. 
Sections \ref{SingularVectors} and \ref{Poisson} explain the relationship between the
residue family operators on spinors and $F$-method
used to produce them in the flat case, resp. the relative Poisson kernel allowing to introduce the integral 
intertwining operator on spinors. In Section \ref{Outlook} we comment on several interesting questions unresolved 
in the article. 
The Appendices A, B and C briefly review the concepts of $Spin$-geometry, Poincar\'e-Einstein metrics and Gegenbauer 
polynomials. 
 
\newpage

\tableofcontents

\section{Boundary value problem for the Dirac operator associated to a Poincar\'e-Einstein metric}\label{Boundary}

  In this section we discuss the construction of formal powers series asymptotic solutions 
  of the eigenvalue equation for the Dirac operator associated to 
  Poincar\'e-Einstein metrics, cf. Appendix B, which we exploit later on in the construction of 
  the residue family operators on spinors. For certain eigenvalues these formal asymptotic solutions are 
  obstructed and the obstructions yield conformal powers of the  Dirac operator, cf. \cite{GMP1, Fischmann}.
  
  Let $(M,h)$ be a $n$-dimensional Riemannian $Spin$-manifold, see Remark \ref{SemiRiemannian} 
  for the discussion of pseudo-Riemannian $Spin$-manifolds. Let $g_+$ be the associated 
  Poincar\'e-Einstein metric on $X=M\times (0,\varepsilon)$, $\varepsilon>0$ and $r$ the 
  coordinate on $[0,\varepsilon)$, cf. Appendix B.
  The conformal compactification of $(X,g_+)$ is 
  \begin{align*}
    (X,\bar{g}:=r^2g_+=dr^2+h_r),
  \end{align*}
  where $\bar{g}$ extends smoothly to $r=0$. Note that we denote 
  by the same letter $X$ the extension of $M\times (0,\varepsilon)$ to $r=0$. 
  The embedding $\iota_r:M\hookrightarrow X$, 
  $M\ni x\mapsto \iota_r(x):=(x,r)$, pulls back 
  $\bar{g}$ to $h_r$, hence $(M,\iota^*_r\bar{g})$ is a hypersurface in $(X,\bar{g})$ with spacelike 
  normal vector field $\partial_r$. Let us denote by 
  \begin{align*}
    S(M,h),\quad S(X,g_+),\quad S(X,\bar{g}),
  \end{align*}
  corresponding spinor bundles, respectively. The Gau\ss{} equation (with respect to $\iota_r$)
  \begin{align*}
    \nabla^{\bar{g}}_Y Z=\nabla^{h_r}_Y Z+\bar{g}(W_r(Y),Z)\partial_r,\quad Y,Z\in\Gamma(TM)
  \end{align*}
  lifts to the spinor bundle
  \begin{align}\label{eq:SpinGauss}
    \nabla_Y^{\bar{g},S}\theta=\widetilde{\nabla}^{h_r}_Y\theta-\frac 12\partial_r\cdot W_r(Y)\cdot\theta,\quad 
    \theta\in\Gamma(S(X,\bar{g})|_{M\times\{r\}}), Y\in\Gamma(TM),
  \end{align}
  where $W_r(Y)=-\nabla^{\bar{g}}_Y\partial_r$ denotes the Weingarten 
  map associated to $\iota_r$. 
  The Dirac operator $\slashed{D}^{\bar{g}}$ on $S(X,\bar{g})$ and the 
  leaf-wise (or, hypersurface) Dirac operator
  \begin{align}
   \widetilde{D}^{h_r}:=\partial_r\cdot\sum_{i=1}^n s_i\cdot\widetilde{\nabla}_{s_i}^{h_r}
     :\Gamma\big(S(X,\bar{g})\big)\to \Gamma\big(S(X,\bar{g})\big)\label{eq:TangDirac}
  \end{align}
  are related by
  \begin{align}
    \iota^*_r\partial_r\cdot\slashed{D}^{\bar{g}}=\widetilde{\slashed{D}}^{h_r}\iota^*_r
      +\frac n2 \iota_r^* H_r-\iota^*_r\nabla^{\bar{g},S}_{\partial_r},\label{eq:HyperDirac}
  \end{align}
  where $H_r=\frac 1n\tr_{h_r}(W_r)$ is the $h_r$-trace of the Weingarten map associated to $\iota_r$
  and $\iota^*_r$, when acting on spinors, denotes the restriction to $M\times \{r\}$, cf. \cite{BGM}. 
  Since $g_+$ and $\bar{g}$ are conformally equivalent metrics, there exists a vector bundle isomorphism
  \begin{align*}
    F_r:S(X,g_+)\to S(X,\bar{g}).
  \end{align*}
  It follows from the conformal covariance of Dirac operator, equation \eqref{eq:HyperDirac} and the 
  isomorphism $F_r$ that $\slashed{D}^{g_+}\varphi=i\lambda\varphi$, $\lambda\in\C$ 
  and $\varphi\in\Gamma\big(S(X,g_+)\big)$, is equivalent to
  \begin{align}
    D(\bar{g})\theta=i\lambda \theta,\quad \theta=F_r\varphi\in\Gamma\big(S(X,\bar{g})\big),   \label{eq:EigenEqn}
  \end{align}
  and 
  \begin{align}
    D(\bar{g}):=-r\partial_r\cdot\widetilde{\slashed{D}}^{h_r}-\frac n2 rH_r \partial_r\cdot
      +r\partial_r\cdot\nabla^{\bar{g},S}_{\partial_r}-\frac n2 \partial_r\cdot\label{eq:DefEigenOp}
  \end{align}
  Let us consider the following formal power series expansion around $r=0$:
  \begin{align*}
    \widetilde{\slashed{D}}^{h_r}=\sum_{k\geq 0}\widetilde{\slashed{D}}^{(h,k)}\frac{r^k}{k!},
      \quad H_r=\sum_{k\geq 0} H^{(k)}\frac{r^k}{k!}.
  \end{align*}
  We shall use the notation $\widetilde{\slashed{D}}:=\widetilde{\slashed{D}}^{h_0}$, see Remark 
  \ref{BundleIdentificationCurved} for its 
  identification with the Dirac operator $\slashed{D}$ on $(M,h)$. 
  The first few terms in the previous expansion depend on $\widetilde{\slashed{D}}$, 
  Schouten tensor $P:=\frac{1}{n-2}(P-Jh)$ and 
  normalized scalar curvature $J:=\frac{\tau}{2(n-1)}$ only, cf. \cite{BGM}:
  \begin{align*}
    \widetilde{\slashed{D}}^{h_r} =\, &\,\tau_0^r\circ\widetilde{\slashed{D}}\circ\tau_r^0
       +\frac 12 r^2\tau_0^r\circ\partial_r\cdot h(P,\widetilde{\nabla}^{h})\circ\tau^0_r+\cdots,\\
    H_r=\, &\frac 1n r J+\frac{1}{3!}r^3\left(-\frac{12}{n}\tr_h(P^2)\right)+\cdots,
  \end{align*}
  where $h(P,\widetilde{\nabla}^{h}):=\sum_{i=1}^n P(s_i)\cdot\widetilde{\nabla}_{s_i}^h$ and 
  $\tau^r_s:\Gamma(S(X,\bar{g})|_{M\times\{s\}})\to \Gamma(S(X,\bar{g})|_{M\times\{r\}})$ denotes the 
  parallel transport along the geodesic coordinate $r$ with respect to $\nabla^{\bar{g},S}$. 
  \begin{bem}\label{BundleIdentificationCurved}
    The restriction of the spin representation $\kappa_{n+1}:Spin(n+1)\to Gl(\Delta_{n+1})$, cf. Appendix A, 
    to $Spin(n)\subset Spin(n+1)$ yields a representation of $Spin(n)$. For even $n$ this is irreducible, 
    whereas for odd $n$ it decomposes into two equivalent irreducible 
    representations. In all cases they are equivalent to $\kappa_n$. Thus there are vector bundle isomorphisms 
    of induced $Spin(n)$-modules: For even $n$, we have 
     \begin{align*}
      \Xi: S(X,\bar{g})|_{r=0}\to S(M,h),
    \end{align*}
    and in case $n$ is odd, we have
    \begin{align*}
      \Xi^\pm: S^\pm(X,\bar{g})|_{r=0}\to S(M,h).
    \end{align*}
    We denoted by $S^\pm(X,\bar{g})$ the splitting of $S(X,\bar{g})$ with respect to the volume element. For 
    any vector field $Y$ on $M$ and spinors $\psi\in S(X,\bar{g})|_{r=0}$ (even $n$)
    or $\psi^\pm\in S^\pm(X,\bar{g})|_{r=0}$ (odd $n$), we have 
    \begin{align*}
       \Xi(\partial_r\cdot Y\cdot\psi)=Y\cdot\Xi(\psi),\quad 
         \Xi^\pm(\pm\partial_r\cdot Y\cdot\psi^\pm)=Y\cdot\Xi^\pm(\psi^\pm),
    \end{align*}
    respectively. In this identification we have $\widetilde{\slashed{D}}=\slashed{D}$ for $n$ even, and 
   $\widetilde{\slashed{D}}=\begin{pmatrix}\slashed{D}&0\\0&-\slashed{D} \end{pmatrix}$ for $n$ odd. Here 
   $\slashed{D}$ is the Dirac operator on $(M,h)$ and by an abuse of notation, we set $\Xi:=\Xi^+$ for $n$ is odd. 
  \end{bem}
  The linear map $\partial_r\cdot:S(X,\bar{g})\to S(X,\bar{g})$ decomposes spinor bundle on $X$ into
  \begin{align*}
    S^{\pm\partial_r}(X,\bar{g}):=\{\theta\in S(X,\bar{g})| \partial_r\cdot\theta=\pm i\theta\},
  \end{align*}
  since it squares to $-1$, and the corresponding projection operators are
  \begin{align}
    \theta^\pm:=P_\pm\theta:=\frac 12(1-(\pm i)\partial_r\cdot)\theta. \label{eq:CurvedProjections}
  \end{align}
  This leads to the definition of convenient function space 
  \begin{align*}
   \mathcal{A}^\pm:=\{\theta=\sum_{j\geq 0} \theta_jr^j| \theta_{2j}\in\Gamma\big(S^{+\partial_r}(X,\bar{g})\big),
      \theta_{2j+1}\in\Gamma\big(S^{-\partial_r}(X,\bar{g})\big),
      \; \nabla^{\bar{g},S}_{\partial_r}\theta_j=0 \},
  \end{align*}
  such that formal asymptotic solutions of equation \eqref{eq:EigenEqn} will be constructed inside $\mathcal{A}^\pm$.  
  \begin{sat}\label{solution-curved}
    Let $\varphi_1\in\Gamma(S^{+\partial_r}(X,\bar{g})|_{r=0})$, 
    $\varphi_2\in\Gamma(S^{-\partial_r}(X,\bar{g})|_{r=0})$ and 
    $\lambda\notin -\N+\frac 12$.
    Then there exists a unique, up to $O(r^{n+1})$ if $n$ is even,
    $\theta\in \mathcal{A}^+$ and $\phi\in \mathcal{A}^-$ such that 
    $r^{\frac n2+\lambda}\theta$ and $r^{\frac n2-\lambda}\phi$ are solutions 
    of equation \eqref{eq:EigenEqn} with $(\theta^+_0)|_{r=0}=\varphi_1$ and 
    $(\phi^-_0)|_{r=0}=\varphi_2$.
  \end{sat}
  \begin{bew}
    Let us consider a formal power series
    \begin{align*}
      \theta=\sum_{j\geq 0} (\theta_j^++\theta^-_j)r^j,
    \end{align*}
    such that $\theta_j^\pm\in\Gamma\big(S^{\pm\partial_r}(X,\bar{g})\big)$, $j\in{\mathbb N}_0$, 
    are parallel along the $r$-coordinate with respect to $\nabla^{\bar{g},S}$. 
    The equation 
    $D(\bar{g})r^{\frac n2+\lambda}\theta=i\lambda r^{\frac n2+\lambda}\theta$ is equivalent to
    \begin{align*}
      i\lambda\sum_{j\geq 0}r^{\frac n2+\lambda+j}(\theta_j^++\theta^-_j)=&
        -\partial_r\cdot\sum_{j,k\geq 0}r^{\frac n2+\lambda+j+k+1}
        \big(\widetilde{\slashed{D}}^{(h,k)}+\frac n2 H^{(k)}\big)(\theta^+_j+\theta^-_j)\\
      &+\partial_r\cdot\sum_{j\geq 0}r^{\frac n2+\lambda+j}(\lambda+j)(\theta^+_j+\theta^-_j).
    \end{align*}
    It follows that $\theta^+_0$ can be chosen arbitrarily, hence we set it to $\varphi_1$ at $r=0$ and 
    extend it over ${X}$ by parallel transport. Furthermore, we have $\theta^-_0=0$. 
    Now the initial data $\varphi_1$ uniquely, up to $O(r^{n+1})$ for even $n$, determine $\theta$ inductively. 

    Similarly, let us consider the formal power series
    \begin{align*}
      r^{\frac n2-\lambda}\phi=\sum_{j\geq 0}r^{\frac n2-\lambda+j}(\phi^+_j+\phi_j^-)
    \end{align*}
    for $\phi_j^\pm\in\Gamma\big(S(X,\bar{g})\big)$, $j\in{\mathbb N}_0$, which are parallel along the $r$-coordinate 
    with respect to $\nabla^{\bar{g},S}$. The initial data  
    $\phi_0^+=0$ and $\phi_0^-$, equal at $r=0$ to $\varphi_2$, determine $\phi$ uniquely 
    up to $O(r^{n+1})$ for even $n$ provided $r^{\frac n2-\lambda}\phi$ solves equation \eqref{eq:EigenEqn}. 
    The proof is complete. 
  \end{bew}

  Let us recall an important application of Proposition \ref{solution-curved}:
  \begin{theorem}\cite[Theorem $2$]{GMP1}

    Let $(M,h)$ be a Riemannian $Spin$-manifold of dimension $n$. For every 
    $N\in\N_0$ ($N<\frac n2$ for even $n$) there exists a conformally covariant linear differential operator
    \begin{align*}
      \mathcal{D}_{2N+1}:\Gamma\big(S(M,h)\big)\to\Gamma\big(S(M,h)\big)
    \end{align*}
    such that $\mathcal{D}_{2N+1}=\slashed{D}^{2N+1}+LOT$, 
    where $LOT$ denotes lower order terms.
  \end{theorem}
  These operators, termed {\sf conformal powers of the Dirac operator}, 
  are the obstruction to extend a boundary spinor into interior, solving 
  equation \eqref{eq:EigenEqn} for $\lambda\in-\N+\frac 12$ recursively.   
  \begin{bem}
    Let us consider the solutions $r^{\frac n2+\lambda}\theta$ and $r^{\frac n2-\lambda}\phi$ 
    of equation \eqref{eq:EigenEqn} in Proposition \ref{solution-curved}.
    By construction, all their coefficients at $r=0$ are
    given by $h$-natural linear differential operators $\mathcal{T}_l(h;\lambda)$, termed {\sf solution operators}, 
    of order $l$, $l\in\N$. They act on boundary data as follows:
    \begin{align*}
       \theta_{2l}^+|_{r=0}=&\mathcal{T}_{2l}(h;\lambda)\varphi_1,\quad  
           \theta_{2l+1}^-|_{r=0}=\,\mathcal{T}_{2l+1}(h;\lambda)\varphi_1
           ,\quad \forall \; l\in\N_0,\\
       \phi_{2l}^-|_{r=0}=&\mathcal{T}_{2l}(h;-\lambda)\varphi_2,\quad
           \phi_{2l+1}^+|_{r=0}=\mathcal{T}_{2l+1}(h;-\lambda)\varphi_2
           ,\quad \forall \; l\in\N_0.
    \end{align*}
    All other coefficients are zero. Note that $\mathcal{T}_l(h;\lambda)$ depends rationally on $\lambda$.
  \end{bem}
  \begin{bsp}
    The first order solution operator is 
    \begin{align*}
        \mathcal{T}_1(h;\lambda)=\frac{1}{2\lambda+1}\widetilde{\slashed{D}},
    \end{align*}
    the second order solution operator is 
    \begin{align*}
        \mathcal{T}_2(h;\lambda)=\frac{1}{2(2\lambda+1)}\widetilde{\slashed{D}}^2
         +\frac 14 J,
    \end{align*}
    and the third order solution operator is 
    \begin{align*}
        \mathcal{T}_3(h;\lambda)
         =&\frac{1}{2(2\lambda+1)(2\lambda+3)}\widetilde{\slashed{D}}^3
         +\frac{1}{2(2\lambda+3)}\partial_r\cdot h(P,\widetilde{\nabla}^{h_r})\\
        &+\frac{1}{4(2\lambda+1)}J\widetilde{\slashed{D}}
         +\frac{1}{4(2\lambda+3)}\partial_r\cdot\grad^M(J)\cdot.
    \end{align*}
    Note that residues of solution operators at $\lambda=-\frac 12,-\frac 32$ yield conformal 
    first and third power of the Dirac operator.
  \end{bsp}  
  
  Now we apply the results obtained so far to the case of flat euclidean manifold 
  $(\R^{n-1},h=\langle\cdot,\cdot\rangle_{n-1})$. 
  The associated Poincar\'e-Einstein metric is given by the 
  hyperbolic metric $g_+={x_n^{-2}}(dx_n^2+h)$ on $\R^n_{x_n>0}$. Its conformal 
  compactification is given by the flat structure $(\R^n_{x_n\geq 0},\bar{g}=\langle\cdot,\cdot\rangle_{n})$,
  and the operator \eqref{eq:TangDirac} simplifies to
  \begin{align}
    \widetilde{\slashed{D}}^h=e_n\cdot\sum_{i=1}^{n-1} e_i\cdot\partial_i.
  \end{align}
  The solution operators $\mathcal{T}_l(h;\lambda)$, $l\in\N_0$, obey an explicit formula
  \begin{align}\label{eq:SolutionOp}
    \mathcal{T}_{2l}(h;\lambda)=\frac{1}{2^{2l} l!(\frac 12+\lambda)_l}\widetilde{\slashed{D}}^{2l},\quad 
      \mathcal{T}_{2l+1}(h;\lambda)=\frac{1}{2^{2l+1} l!(\frac 12+\lambda)_{l+1}}\widetilde{\slashed{D}}^{2l+1},
  \end{align}
  where $(a)_l:=a(a+1)\ldots(a+l-1)$, $a\in\C$ and $l\in\N$, denotes the Pochhammer symbol. 
  We set conventionally $(a)_0:=1$.    
  
  \begin{bem}(Pseudo-Riemannian structures)\label{SemiRiemannian}

     Let $(M^{n},h)$ be a pseudo-Riemannian $Spin$-manifold of signature $(p,q)$, $p+q=n$. 
     In this case we can decompose $S(X,\bar{g})$ into 
     $(\pm i)^{p+1}$-eigenspaces with respect to $\partial_r$. The projector operators are defined by
     \begin{align*}
        P_\pm:=\frac 12(1-(\pm i)^{p+1}\partial_r\cdot),
     \end{align*}
     hence formally self-adjoint with respect to the spinor scalar product $<\cdot,\cdot>$ on $S(X,\bar{g})$.
     This is compatible with the eigenvalue equation analogous to
     \eqref{eq:EigenEqn}, replacing $i\lambda$ by $i^{p+1}\lambda$, which does
     not change the structure of formal solutions obtained in Proposition \ref{solution-curved}. 
  \end{bem}


\section{Residue family operators on spinors}\label{ResidueFamilies}

The present section introduces a definition of residue family operators on spinors. 
Firstly, we define the residue family operators on spinors for general semi-Riemannian 
$Spin$-manifolds. Secondly, we prove that in the flat case the resulting residue family operators 
agree with a family of intertwining operators introduced in \cite{KOSS}. 

\subsection{Residue family operators on spinors - curved case}
  The present subsection is devoted to the general definition of the residue family operators on spinors 
  for semi-Riemannian $Spin$-manifolds. Again, we restrict ourselves to Riemannian $Spin$-manifolds and 
  discuss the pseudo-Riemannian case in Remark \ref{PseudoRiemannianCase}. 

  Let us introduce the volume function 
  \begin{align*}
    v(r):=\sqrt{\frac{\det(h_r)}{\det(h)}}=1+r^2v_2+r^4v_4+\cdots,
  \end{align*}
  where $v_{2j}$, $j>0$, are called renormalized volume coefficients, see \cite{Graham}.
  We have the following relation among volume forms,
  \begin{align*}
    Vol(\bar{g})=v(r)drVol(h).
  \end{align*}
  Now, we consider the formal asymptotic solution 
  \begin{align*}
    \bar{\theta}:=r^{\frac n2+\lambda}\theta=\sum_{j\geq 0} r^{\frac n2+\lambda+j}\theta_j
  \end{align*}
  of the eigenvalue equation \eqref{eq:EigenEqn} for 
  the eigenvalue $i\lambda$, cf. Proposition \ref{solution-curved}. 
  For $\mu\in\C$ such that $Re(\mu)\gg 0$, we define a family of spinor valued distributions by
  \begin{align*}
    M_{\bar{\theta}}(\mu;r)(\varphi):=\int_{[0,\varepsilon)}\int_M <r^\mu\bar{\theta},\varphi>Vol(\bar{g}),
  \end{align*}
  where $\varphi\in\Gamma_c\big(S(X,\bar{g})\big)$ is compactly supported 
  near $r=0$ and $<\cdot,\cdot>$ denotes the 
  scalar product on the spinor bundle $S(X,\bar{g})$.

  After fixing the eigenvalue $\lambda$, we can meromorphically extend $M_{\bar{\theta}}(\mu;r)$ to $\mu\in\C$ 
  with simple poles at $\mu\in -\frac n2-\lambda-1-\N_0$. 
  By partial integration with respect to $\nabla^{\bar{g},S}_{\partial_r}$, 
  the residue of $M_{\bar{\theta}}(\mu;r)$ at $\mu=-\frac n2-\lambda-1-N$, $N\in\N_0$ 
  ($N\leq n$ for even $n$), is
  \begin{align*}
    &Res_{\mu=-\frac n2-\lambda-1-N}\big(M_{\bar{\theta}}(\mu;r)(\varphi)\big)=\\
    &\int_M\sum_{j=0}^N\frac{1}{(N-j)!}<\sum_{k=0}^jv_k\mathcal{T}_{j-k}(h;\lambda)\theta_0^+
     ,\big((\nabla^{\bar{g},S}_{\partial_r})^{N-j}\varphi\big)(0,\cdot)>Vol(h).
  \end{align*}
  Taking adjoints (with respect to the induced $L^2$-scalar product), 
  indicated by $^\star$, of the operators on the left hand side of the scalar product $<\cdot,\cdot>$, 
  we arrive at a family of linear $h$-natural differential operators  
  \begin{align*}
    \slashed{\delta}^+_N(h;\lambda):\Gamma_c\big(S(X,\bar{g})\big)\to \Gamma(S^{+\partial_r}(X,\bar{g})|_{r=0})
  \end{align*}
  defined by 
  \begin{align*}
    Res_{\mu=-\frac n2-\lambda-1-N}\big(M_{\bar{\theta}}(\mu;r)(\varphi)\big)
     =\int_M<\theta_0^+,\slashed{\delta}^+_N(h;\lambda)\varphi>Vol(h).
  \end{align*}
  They are given by explicit formulas 
  \begin{align}
    \slashed{\delta}^+_{2N}(h;\lambda)
      =&\sum_{j=0}^{N-1}\frac{1}{(2N-2j-1)!}
        \sum_{k=0}^j\big(\mathcal{T}_{2j+1-2k}(h;\lambda)\big)^\star 
         v_{2k}\iota^*(\nabla^{\bar{g},S}_{\partial_r})^{2N-2j-1}\notag\\
      &+\sum_{j=0}^N\frac{1}{(2N-2j)!}
        \sum_{k=0}^j\big(\mathcal{T}_{2j-2k}(h;\lambda)\big)^\star 
        v_{2k}\iota^*(\nabla^{\bar{g},S}_{\partial_r})^{2N-2j},\label{eq:EvenDistr}\\
    \slashed{\delta}^+_{2N+1}(h;\lambda)
      =&\sum_{j=0}^N\frac{1}{(2N-2j+1)!}
        \sum_{k=0}^j\big(\mathcal{T}_{2j-2k}(h;\lambda)\big)^\star
        v_{2k}\iota^*(\nabla^{\bar{g},S}_{\partial_r})^{2N+1-2j}\notag\\
      &+\sum_{j=0}^{N}\frac{1}{(2N-2j)!}
        \sum_{k=0}^j\big(\mathcal{T}_{2j+1-2k}(h;\lambda)\big)^\star 
        v_{2k}\iota^*(\nabla^{\bar{g},S}_{\partial_r})^{2N-2j},\label{eq:OddDistr}
  \end{align}
  where $\iota:=\iota_0$ denotes the embedding of $M$ into $M\times [0,\varepsilon)$ corresponding to $r=0$. 
  Similarly, one defines for $N\in{\mathbb N}_0$ 
  \begin{align*}
     \slashed{\delta}^-_N(h;\lambda):\Gamma_c\big(S(X,\bar{g})\big)\to \Gamma(S^{-\partial_r}(X,\bar{g})|_{r=0})
  \end{align*}
  by 
  \begin{align*}
    Res_{\mu=-\frac n2-\lambda-1-N}\big(M_{\bar{\phi}}(\mu;r)(\varphi)\big)
     =\int_M<\phi_0^-,\slashed{\delta}^-_N(h;\lambda)\varphi>Vol(h),
  \end{align*}
  where $\bar{\phi}:=r^{\frac n2+\lambda}\phi$ is a solution of the eigenvalue equation 
  \eqref{eq:EigenEqn} for the eigenvalue $-i\lambda$. Again, 
  $\slashed{\delta}^\pm_N(h;\lambda)$ are structurally the same except 
  that their target spaces are different. We define, by continuous extension $\slashed{\delta}^\pm_N(h;\lambda)$ 
  to the space of smooth sections,  
  \begin{align}
    \slashed{\delta}_N(h;\lambda):=\slashed{\delta}^+_N(h;\lambda)+\slashed{\delta}^-_N(h;\lambda)
      :\Gamma\big(S(X,\bar{g})\big)\to \Gamma(S(X,\bar{g})|_{r=0}).
  \end{align}
  It follows from equations \eqref{eq:EvenDistr} and \eqref{eq:OddDistr} that 
  $\slashed{\delta}_{2N}(h;\lambda)$ has $N$ simple poles at $\lambda\in\{-N+\frac 12,\ldots,-\frac 12\}$, whereas 
  $\slashed{\delta}_{2N+1}(h;\lambda)$ has $(N+1)$ simple poles at $\lambda\in\{-N-\frac 12,\ldots,-\frac 12\}$. 

  Now we define the residue family operators for spinors.
  \begin{defn}\label{DefResFamOp}
    Let $\lambda\in \C$ and $N\in\N_0$ ($N\leq n$ for even $n$). Then the differential operators 
    \begin{align*}
      \slashed{D}^{res}_{N}(h;\lambda):\Gamma\big(S(X,\bar{g})\big)\to\Gamma(S(X,\bar{g})|_{r=0})
    \end{align*}
    given by 
    \begin{align*}
      \slashed{D}^{res}_{2N}(h;\lambda):=&2^{2N}N!\left[(\lambda+\frac n2-2N+\frac 12)\cdot \ldots 
       \cdot(\lambda+\frac n2-N-\frac 12)\right]\times\\
      &\times\slashed{\delta}_{2N}(h;\lambda+\frac{n}{2}-2N),\\
      \slashed{D}^{res}_{2N+1}(h;\lambda):=&2^{2N+1}N!\left[(\lambda+\frac n2-2N-\frac 12)\cdot\ldots
        \cdot(\lambda+\frac n2-N-\frac 12)\right]\times\\
      &\times\slashed{\delta}_{2N+1}(h;\lambda+\frac n2-2N-1),
    \end{align*}
    are called the residue family operators on spinors.
  \end{defn}
  All statements about residue family operators $\slashed{D}^{res}_{N}(h;\lambda)$ for even $n$ are in general
  valid only for $N\leq \frac n2$. In what follows this restriction will be omitted. 

  \begin{bem}
    The last definition is analogous to the definition of residue families in 
    the scalar case, cf. \cite{Juhl}. However,
    contrary to the scalar case, the target space of residue family operators is $S(X,\bar{g})|_{r=0}$ 
    rather than $S(M,h)$. An explanation for this choice will be given in Remark \ref{WhyNoIdentification}. 
  \end{bem}
  \begin{bem}
    It follows from the definition of residue family operators:
    \begin{enumerate}
      \item $\slashed{D}_N^{res}(h;\lambda)$ is a family of linear $h$-natural differential operators 
       of order $N$.
      \item $\slashed{D}_N^{res}(h;\lambda)$ is a polynomial in 
      $\lambda$ of degree $\lfloor\frac N2\rfloor+1$, where $\lfloor k\rfloor$ denotes the integer part of $k\in\R_+$. 
    \end{enumerate} 
  \end{bem}
  We present a few low order examples.
  \begin{bsp}\label{ExamplesCurved}
    The first order residue family operator is 
    \begin{align*}
       \slashed{D}_1^{res}(h;\lambda)=&\, 2(\lambda+\frac n2-\frac 12)\slashed{\delta}_1(h;\lambda+\frac n2-1)\\
       =&\, (2\lambda+n-1)[\iota^*\nabla^{\bar{g},S}_{\partial_r}
          +\big(\mathcal{T}_1(h;\lambda+\frac n2-1)\big)^\star\iota^*].
    \end{align*}
    Because the adjoint of solution operator is 
    \begin{align*}
      \big(\mathcal{T}_1(h;\lambda)\big)^\star=\frac{1}{2\lambda+1}\widetilde{\slashed{D}},
    \end{align*} 
    we get 
    \begin{align*}
      \slashed{D}_1^{res}(h;\lambda)
         =&\, (2\lambda+n-1)\iota^*\nabla^{\bar{g},S}_{\partial_r}+\widetilde{\slashed{D}}\iota^*.
    \end{align*}
    The second order residue family operator is 
    \begin{align*}
       \slashed{D}_2^{res}(h;\lambda)=&\, 2^2(\lambda+\frac n2-\frac 32)
          \slashed{\delta}_2(h;\lambda+\frac n2-2)\\
       =&\, 2(2\lambda+n-3)\big[\frac 12\iota^*(\nabla^{\bar{g},S}_{\partial_r})^2
          +\big(\mathcal{T}_2(h;\lambda+\frac n2-2)\big)^\star\iota^*\\
        &\, +v_2\iota^*+\big(\mathcal{T}_1(h;\lambda+\frac n2-2)\big)^\star\iota^*
          \nabla^{\bar{g},S}_{\partial_r}\big].
    \end{align*}
    Using $v_2=-\frac 12 J$, cf. \cite[Theorem $6.9.2$]{Juhl}, and 
    \begin{align*}
       \big(\mathcal{T}_2(h;\lambda)\big)^\star=\frac{1}{2(2\lambda+1)}\widetilde{\slashed{D}}^2+\frac 14 J,
    \end{align*}
    we have 
    \begin{align*}
      \slashed{D}_2^{res}(h;\lambda)=&\, (2\lambda+n-3)\iota^*(\nabla^{\bar{g},S}_{\partial_r})^2
        +\widetilde{\slashed{D}}^2\iota^*
         -\frac{1}{2}({2\lambda+n-3})J\iota^*+2\widetilde{\slashed{D}}\iota^*\nabla^{\bar{g},S}_{\partial_r}.
    \end{align*}
    The third order residue family operator is 
    \begin{align*}
       \slashed{D}_3^{res}(h;\lambda)=&\, 2^3(\lambda+\frac n2-\frac 52)(\lambda+\frac n2-\frac 32)
           \slashed{\delta}_3(h;\lambda+\frac n2-3)\\
       =&\, 2(2\lambda+n-5)(2\lambda+n-3)\big[\frac{1}{3!}\iota^*(\nabla^{\bar{g},S}_{\partial_r})^3
          +\big(\mathcal{T}_2(h;\lambda+\frac n2-3)\big)^\star\iota^*\nabla^{\bar{g},S}_{\partial_r}\\
        &\, +v_2\iota^*\nabla^{\bar{g},S}_{\partial_r}
           +\frac 12 \big(\mathcal{T}_1(h;\lambda+\frac n2-3)\big)^\star\iota^*(\nabla^{\bar{g},S}_{\partial_r})^2\\
        &\, +\big(\mathcal{T}_3(h;\lambda+\frac n2-3)\big)^\star\iota^*
           +\big(\mathcal{T}_1(h;\lambda+\frac n2-3)\big)^\star v_2\iota^*\big].
    \end{align*}
    Using  
    \begin{align*}
       \big(\mathcal{T}_3(h;\lambda)\big)^\star=&\, \frac{1}{2(2\lambda+1)(2\lambda+n-3)}
         \widetilde{\slashed{D}}^3
         +\frac{1}{2(2\lambda+n-3)}\partial_r\cdot h(P,\widetilde{\nabla}^h)\\
       &\, +\frac{1}{2(2\lambda+n-3)}\partial_r\cdot\grad^M(J)\cdot
         +\frac{1}{4(2\lambda+1)}J\widetilde{\slashed{D}}\\
       &\, +\frac{1}{2(2\lambda+1)(2\lambda+n-3)}\partial_r\cdot\grad^M(J)\cdot,
    \end{align*}
    where we note that $\partial_r\cdot h(P,\widetilde{\nabla}^h)$ is not self-adjoint because its adjoint contributes 
    by the gradient of $J$, we get
    \begin{align*}
      \slashed{D}_3^{res}(h;\lambda)
        =&\frac{1}{3}{(2\lambda+n-5)(2\lambda+n-3)}\iota^*(\nabla^{\bar{g},S}_{\partial_r})^3
        +(2\lambda+n-3)\widetilde{\slashed{D}}^2\iota^*\nabla^{\bar{g},S}_{\partial_r}\\
      &-\frac{1}{2}{(2\lambda+n-5)(2\lambda+n-3)}J\iota^*\nabla^{\bar{g},S}_{\partial_r}
        +(2\lambda+n-3)\widetilde{\slashed{D}}\iota^*(\nabla^{\bar{g},S}_{\partial_r})^2\\
      &+\widetilde{\slashed{D}}^3\iota^*+(2\lambda+n-5)\partial_r\cdot h(P,\widetilde{\nabla}^h)\iota^*\\
      &-\frac{1}{2}({2\lambda+n-3})J\widetilde{\slashed{D}}\iota^*-\partial_r\cdot\grad^M(J)\cdot\iota^*.
    \end{align*}
    Later on these formulas will be used to prove the factorization identities in the curved case. 
  \end{bsp}

  \begin{bem}\label{PseudoRiemannianCase}
    Let $(M,h)$ be a pseudo-Riemannian $Spin$-manifold. In Remark \ref{SemiRiemannian} we suggested how to 
    solve the eigenvalue equation in the presence of signature. The adjustments made there can be used to define 
    residue family operators in a similar way as for the Riemannian $Spin$-manifolds. 
  \end{bem}


\subsection{Residue family operators on spinors - flat case}\label{ResFamFlatCase}
  We now specialize the definition of residue family operators $\slashed{D}^{res}_N(h;\lambda)$, 
  cf. Definition \ref{DefResFamOp}, to the flat case $(\R^{n-1},h)$. 
  Notice that in this case there is no restriction on $N\in\N_0$ for even $(n-1)$. Due to explicit knowledge 
  of solution operators \eqref{eq:SolutionOp}, we are able to derive an explicit formula for residue 
  family operators in terms of Gegenbauer polynomials.  

  First note that formal self-adjointness 
  (with respect to the $L^2$-scalar product $<\cdot,\cdot>_{L^2}$ on $S(\R^n,\bar{g})$) of the operator  
  \begin{align*}
   \widetilde{\slashed{D}}=e_n\cdot \sum_{i=1}^{n-1} e_i\cdot\partial_i,
  \end{align*} 
  i.e., $<\widetilde{\slashed{D}}\theta_1,\theta_2>_{L^2}=<\theta_1,\widetilde{\slashed{D}}\theta_2>_{L^2}$ 
  for compactly supported spinors $\theta_1,\theta_2$, implies that solution operators, cf. equation 
  \eqref{eq:SolutionOp}, are formally self-adjoint, 
  $\mathcal{T}_{l}(h;\lambda)^\star=\mathcal{T}_{l}(h;\lambda)$ for $l\in\N_0$. 
  Secondly, we introduce 
  \begin{align*}
    D_{\mathcal{T}}&:=\sum_{i=1}^{n-1}e_i\cdot\partial_i,\quad\text{(tangential Dirac operator)},\\
    D_{\mathcal{N}}&:=e_n\cdot\partial_n,\quad  \text{(normal Dirac operator)}.
  \end{align*} 
  The even and odd order family \eqref{eq:EvenDistr} and \eqref{eq:OddDistr} specialize, 
  using $v(r)=1$ and $\nabla^{\bar{g},S}_{\partial_r}=\partial_n$ in the flat case, to
  \begin{align}
   \slashed{\delta}&_{2N}^+(h;\lambda)
     =\sum_{j=0}^N\frac{1}{(2N-2j)!}
       \left(\frac{1}{4^j j!(\frac 12+\lambda)_j}\widetilde{\slashed{D}}^{2j}\right)
       \iota^*\partial_n^{2N-2j}\notag\\
     &-\sum_{j=0}^{N-1}\frac{1}{(2N-2j-1)!}
      \left(\frac{1}{2^{2j+1}j!(\frac 12+\lambda)_{j+1}}\widetilde{\slashed{D}}^{2j}\right)
      \iota^*\partial_n^{2N-2j-2}D_{\mathcal{T}}D_{\mathcal{N}},\label{eq:EvenFlatDistr}\\
    \slashed{\delta}&_{2N+1}^+(h;\lambda)
     =e_n\cdot\bigg[\sum_{j=0}^{N}\frac{1}{(2N-2j)!}
      \left(\frac{1}{2^{2j+1} j! (\frac 12+\lambda)_{j+1}}\widetilde{\slashed{D}}^{2j}\right)
      \iota^*\partial_n^{2N-2j}D_{\mathcal{T}}\notag\\
     &-\sum_{j=0}^N\frac{1}{(2N-2j+1)!}
       \left(\frac{1}{4^j j! (\frac 12+\lambda)_j}\widetilde{\slashed{D}}^{2j}\right)
       \iota^*\partial_n^{2N-2j}D_{\mathcal{N}}\bigg].\label{eq:OddFlatDistr}
  \end{align}
  where $\iota:\R^{n-1}\to \R^n$ is canonical embedding on first $n-1$ coordinates. 
  Note that similar statement holds for $\delta^-_N(h;\lambda)$. 
 
  For $N\in\N_0$, $\lambda\in\C$ and $a_N^{(N)}(\lambda),b_N^{(N)}(\lambda)\in\R$, we define 
  \begin{align}
    a^{(N)}_j(\lambda):=&\frac{N! (-2)^{N-j}}{j!(2N-2j)!}
      \prod_{k=j}^{N-1}(2\lambda-4N+2k+n+1)a_N^{(N)}(\lambda),\label{eq:DefA}\\
    b^{(N)}_j(\lambda):=&\frac{N! (-2)^{N-j}}{j!(2N-2j+1)!}
      \prod_{k=j}^{N-1}(2\lambda-4N+2k+n-1)b_N^{(N)}(\lambda),\label{eq:DefB}
  \end{align}
  $0\leq j\leq N-1$. They satisfy the following recurrence relations:
  \begin{align}
    (N-j+1)(2N-2j+1)a^{(N)}_{j-1}(\lambda)+j(2\lambda+n-4N+2j-1)a_j^{(N)}(\lambda)&=0,\label{eq:a-recurrence}\\
    (N-j+1)(2n-2j+3)b_{j-1}^{(N)}(\lambda)+j(2\lambda+n-4N+2j-3)b_j^{(N)}(\lambda)&=0,\label{eq:b-recurrence}
  \end{align}
   for all $1\leq j\leq N$ and $\lambda\in\C$ and are known as the coefficients of 
  Gegenbauer polynomials of even and odd 
  degree, respectively; cf. Appendix C for the origin and wider framework of this notion. 
  We use the conventions $a^{(N)}_{-1}:=0$ and $b^{(N)}_{-1}:=0$.
  Note that empty products are set to be $1$, and we shall consider the 
  normalizations $a_N^{(N)}(\lambda)=b_N^{(N)}(\lambda)=(-1)^N$. 
  \begin{bem}
    The normalization among odd and even Gegenbauer polynomials is a 
    consequence of equations \eqref{eq:EvenNorm1}, \eqref{eq:EvenNorm2} 
    and \eqref{eq:OddNorm1}, \eqref{eq:OddNorm2}.    
  \end{bem}
  The next Theorem identifies $\slashed{D}_N^{res}(h;\lambda)$ 
  as a family of intertwining differential operators $D_N(\lambda)$ introduced in \cite{KOSS}, though
  we use the opposite convention in the sign of the density given by $\lambda$.
  \begin{theorem}\label{Intertwiner}
    Let $\lambda\in\C$ and $N\in\N_0$. Then it holds
    \begin{align}
      \slashed{D}^{res}_{2N}(h;\lambda)=(-1)^N\Big[&\sum_{j=0}^Na_j^{(N)}
        (\lambda-\frac 12)D_{\mathcal{T}}^{2j}\iota^*
        D_{\mathcal{N}}^{2N-2j}\notag \\
      &+2N\sum_{j=0}^{N-1}b_j^{(N-1)}
        (\lambda-\frac 12)D_{\mathcal{T}}^{2j+1}\iota^*D_{\mathcal{N}}^{2N-2j-1}\Big],\\
      \slashed{D}^{res}_{2N+1}(h;\lambda)=(-1)^{N+1}e_n\cdot
        \Big[&\sum_{j=0}^N b_j^{(N)}(\lambda-\frac 12)c^{(N)}(\lambda)
        D_{\mathcal{T}}^{2j}\iota^*D_{\mathcal{N}}^{2N-2j+1}\notag \\
      &-\sum_{j=0}^N a_j^{(N)}(\lambda-\frac 12)D_{\mathcal{T}}^{2j+1}\iota^*D_{\mathcal{N}}^{2N-2j} \Big],
    \end{align}
    where $c^{(N)}(\lambda):=(2\lambda+n-2N-2)$.
  \end{theorem}
  \begin{bew}
   The proof is based on direct computations.
   Let us start with the even order family. From equation \eqref{eq:EvenFlatDistr} we obtain
   \begin{align*}
    \slashed{D}^{res}_{2N}&(h;\lambda)=2^{2N}N!\left[(\lambda+\frac n2-2N)\cdot \ldots 
       \cdot(\lambda+\frac n2-N-1)\right]\times\\
      &\times\slashed{\delta}_{2N}(h;\lambda+\frac n2-2N-\frac 12)\\
    =&2^{2N}N!\left[(\lambda+\frac n2-2N)\cdot \ldots 
       \cdot(\lambda+\frac n2-N-1)\right]\times\\
      &\times\Big[\sum_{j=0}^N\frac{1}{(2N-2j)!}\frac{1}{2^{2j}j!(\lambda+\frac n2-2N)_j}
         D_{\mathcal{T}}^{2j}\iota^*\partial_n^{2N-2j}\\
    &\quad -\sum_{j=0}^{N-1}\frac{1}{(2N-2j-1)!}\frac{1}{2^{2j+1}j!(\lambda+\frac n2-2N)_{j+1}}
         D_{\mathcal{T}}^{2j}\iota^*\partial_n^{2N-2j-2}D_{\mathcal{T}}D_{\mathcal{N}}  \Big]\\
    =&\sum_{j=0}^N\frac{2^{2N-2j}N!}{(2N-2j)!j!}
         \prod_{k=j}^{N-1}(\lambda+\frac n2-2N+k)D_{\mathcal{T}}^{2j}\iota^*\partial_n^{2N-2j}\\
    &-\sum_{j=0}^{N-1}\frac{2^{2N-2j-1}N!}{(2N-2j-1)!j!}
       \prod_{k=j}^{N-2}(\lambda+\frac n2-2N+k+1)
       D_{\mathcal{T}}^{2j}\iota^*\partial_n^{2N-2j-2}D_{\mathcal{T}}D_{\mathcal{N}}.
   \end{align*}
   Introducing $1=(-1)^{l}e_n^{2l}$ for $l=N-j\in\N$ and using definitions \eqref{eq:DefA} 
   and \eqref{eq:DefB}, we get
   \begin{align*}
    \slashed{D}^{res}_{2N}(h;\lambda)
     =(-1)^N\Big[&\sum_{j=0}^Na_j^{(N)}(\lambda-\frac 12)D_{\mathcal{T}}^{2j}\iota^*
        D_{\mathcal{N}}^{2N-2j}\\
     &+2N\sum_{j=0}^{N-1}b_j^{(N-1)}
      (\lambda-\frac 12)D_{\mathcal{T}}^{2j+1}\iota^*D_{\mathcal{N}}^{2N-2j-1}\Big].
   \end{align*}
   The odd order family is treated in an analogous way and we obtain
   \begin{align*}
    \slashed{D}^{res}_{2N+1}(h;\lambda)
    =(-1)^{N+1}e_n\cdot\Big[&\sum_{j=0}^N b_j^{(N)}(\lambda-\frac 12)(2\lambda+n-2N-2)
     D_{\mathcal{T}}^{2j}\iota^*D_{\mathcal{N}}^{2N-2j+1}\\
    &\quad-\sum_{j=0}^N a_j^{(N)}(\lambda-\frac 12)D_{\mathcal{T}}^{2j+1}\iota^*D_{\mathcal{N}}^{2N-2j} \Big],
   \end{align*}
   which completes the proof.
  \end{bew}

  \begin{bsp}\label{ExamplesFlat}
   Let us give several low order examples:
   \begin{align*}
    \slashed{D}^{res}_{1}(h;\lambda)
    =&e_n\cdot\big[D_{\mathcal{T}}\iota^*-(2\lambda+n-2)\iota^*D_{\mathcal{N}}\big],\\
    \slashed{D}^{res}_{2}(h;\lambda)
    =&D_{\mathcal{T}}^2\iota^*-2D_{\mathcal{T}}\iota^*D_{\mathcal{N}}
      -(2\lambda+n-4)\iota^*D_{\mathcal{N}}^2,\\
    \slashed{D}^{res}_{3}(h;\lambda)
    =&e_n\cdot\big[D_{\mathcal{T}}^3\iota^*-(2\lambda+n-4)(D_{\mathcal{T}}^2\iota^*D_{\mathcal{N}}
      +D_{\mathcal{T}}\iota^*D_{\mathcal{N}}^2)\\
    &+\frac 13 (2\lambda+n-6)(2\lambda+n-4)\iota^*D_{\mathcal{N}}^3  \big],\\
    \slashed{D}^{res}_{4}(h;\lambda)
    =&D_{\mathcal{T}}^4\iota^*-4 D_{\mathcal{T}}^3\iota^*D_{\mathcal{N}}
       +\frac 43 (2\lambda+n-6)D_{\mathcal{T}}\iota^*D_{\mathcal{N}}^3\\
    &-2(2\lambda+n-6)D_{\mathcal{T}}^2\iota^*D_{\mathcal{N}}^2
      +\frac 13(2\lambda+n-8)(2\lambda+n-6)\iota^*D_{\mathcal{N}}^4.
   \end{align*}
  \end{bsp}


\section{Conformal transformation properties}\label{ConformalTrafo}

  An important property of the residue family operators on spinors $\slashed{D}^{res}_{N}(h;\lambda)$,
  implicit in the invariance properties of solution operators in the original eigenvalue 
  problem for the Poincar\'e-Einstein Dirac operator, is its conformal covariance with respect to a conformal
  change of the boundary metric. The conformal covariance property determines the 
  residue family of natural differential operators on spinors in the flat case uniquely.   

  We shall start with several remarks used to understand the conformal covariance property of the residue
  family operators. Let $(X_1,g_1)$, $(X_2,g_2)$ be semi-Riemannian $Spin$-manifolds of dimension $n$, and 
  $\Phi: X_1\to X_2$
  a $Spin$-structure preserving diffeomorphism such that $\Phi^*(g_2)=g_1$. 
  We denote by $\Phi^*$ the pull-back and by 
  $\Phi_*$ the push-forward map induced by $\Phi$. For example, $\Phi$ pulls-back sections of associated 
  vector bundles on $X_2$ to $X_1$. 
  Note that the diffeomorphism $\Phi$ lifts to a vector bundle isomorphism 
  \begin{align*}
     \widetilde{\Phi}:S(X_1,g_1)\to S(X_2,g_2),
  \end{align*}
  and it is compatible with Clifford multiplication 
  $\widetilde{\Phi}(X\cdot \psi)=\Phi_* X\cdot \widetilde{\Phi}\psi$, $X\in TX_1$ and $\psi\in S(X_1,g_1)$. 
  For $\varphi\in \Gamma\big(S(X_2,g_2)\big)$
  we set $\widetilde{\Phi}^*\varphi:=\widetilde{\Phi}^{-1}\circ \varphi\circ\Phi\in \Gamma\big(S(X_1,g_1)\big)$, 
  which is the pull-back of $\varphi$ with respect to $\widetilde{\Phi}$. 
 
  \begin{lem}
    In the setting described above, the Levi-Civita connections are related by 
    \begin{align}
       \nabla^{g_1}=&(\Phi^{-1})_*\nabla^{g_2}\Phi_*,\label{eq:CovDiff}
    \end{align}
    and induced spinor covariant derivatives are related by
    \begin{align}
       (\widetilde{\Phi}^{-1})^*\nabla_X^{g_1,S}\widetilde{\Phi}^*\varphi
          =&\nabla^{g_2,S}_{\Phi_* X}\varphi,\label{diffd}
    \end{align}
    for $\varphi\in \Gamma\big(S(X_2,g_2)\big)$, $X\in TX_1$. 

    As for the Dirac operators on $(X_1,g_1)$ and $(X_2,g_2)$, we have 
    \begin{align}
       (\widetilde{\Phi}^{-1})^*\slashed{D}^{g_1}\widetilde{\Phi}^*\varphi
          =\slashed{D}^{g_2}\varphi,\label{eq:diffDirac}
    \end{align}
    for $\varphi\in\Gamma\big(S(X_2,g_2)\big)$.
  \end{lem}
  \begin{bew}
    The first statement is the isometry invariance of Levi-Civita connections. 
    The second claim follows from the local description of spinor 
    covariant derivatives, equation \eqref{eq:CovDiff} 
    and $(\tilde{\Phi}^{-1})^*(X\cdot \tilde{\Phi}^*\varphi)=\Phi_*X\cdot\varphi$ for $X\in\Gamma(TM)$ and 
    $\varphi\in\Gamma\big(S(X_2,g_2)\big)$. Finally, using equation \eqref{diffd} we get for Dirac operators
    \begin{align*}
      (\tilde{\Phi}^{-1})^*\slashed{D}^{g_1}\tilde{\Phi}^*\varphi
        =&(\tilde{\Phi}^{-1})^*\sum_{i=1}^n\varepsilon_i s_i\cdot \nabla^{g_1,S}_{s_i}\tilde{\Phi}^*\varphi 
      =(\tilde{\Phi}^{-1})^*
        \sum_{i=1}^n\varepsilon_i s_i\cdot\tilde{\Phi}^* \nabla^{g_2,S}_{\Phi_* s_i}\varphi\\
      =&\sum_{i=1}^n\varepsilon_i \Phi_* s_i\cdot
        \nabla^{g_2,S}_{\Phi_* s_i}\varphi=\slashed{D}^{g_2}\varphi,
    \end{align*}
    which completes the proof. 
  \end{bew}
  
  These results will be applied to the diffeomorphism of the Poincar\'e-Einstein metric 
  induced by a conformal change of the boundary metric $h$, cf. Appendix B.

  Let $g_+=r^{-2}(dr^2+h_r)$ be the Poincar\'e-Einstein metric on $X$ for a representative 
  $h\in[h]$ in the conformal class on the conformal infinity $(M,[h])$ of $X$. 
  For a smooth function $\sigma\in C^\infty(M)$, the Poincar\'e-Einstein
  metrics of $h$ and $\widehat{h}=e^{2\sigma}h$ are related by 
  \begin{align*}
    \Phi^*\big(r^{-2}(dr^2+h_r)\big)=r^{-2}(dr^2+\widehat{h}_r),
  \end{align*}
  where $\Phi$ is the induced diffemorphism on $X$
  which restricts to the identity on the hypersurface $r=0$, cf. Appendix B. 
  Then we have 
  \begin{align}\label{pemetrictrans}
     \Phi^*(dr^2+h_r)=\left(\frac{\Phi^*(r)}{r}\right)^2(dr^2+\widehat{h}_r)
  \end{align}
  with 
  \begin{align*}
    \lim_{r\to 0}\frac{\Phi^*(r)}{r}=e^{-\sigma},
      \quad \iota^*\left(\frac{\Phi^*(r)}{r}\right)^{\mu}=e^{-\mu\sigma},
  \end{align*}
  where $\iota$ is the usual embedding of $M$ into $X$ at $r=0$. 
  \begin{theorem}\label{ConformalTrafoLawResidueFamily}
    The residue family operators $\slashed{D}^{res}_{N}(h;\lambda)$ are 
    conformally covariant differential operators in the sense that
    \begin{align*}
       \slashed{D}^{res}_{N}(\widehat{h};\lambda)=e^{(\lambda-N)\varphi}\circ
         \slashed{D}^{res}_{N}({h};\lambda)\circ(\tilde{\Phi}^{-1})^*\circ\left(\frac{\Phi^*(r)}{r}\right)^\lambda 
    \end{align*} 
    for all $\sigma\in C^\infty(M)$, $N\in{\mathbb N}_0$.
  \end{theorem}
  \begin{bem}
    Notice that the behavior of $\slashed{D}^{res}_{N}({h};\lambda)$ under the conformal transformation is actually not 
    conformally covariant in the usual sense, cf. \cite{Kosmann}. 
  \end{bem}
  \begin{bew}
    Let us consider the conformal compactification $(M\times [0,\varepsilon),\bar{g})$ 
    of the Poincar\'e-Einstein metric $g_+$ associated to $h\in[h]$ on $M$, 
    $\bar{\theta}:=r^{\frac n2+\lambda}\theta\in\Gamma\big(S(X,\bar{g})\big)$ a formal asymptotic 
    solution of the eigenvalue 
    equation \eqref{eq:EigenEqn}, cf. Proposition \ref{solution-curved}, 
    and $\phi\in\Gamma_c\big(S(X,\bar{g})\big)$ a smooth 
    compactly supported spinor field. 

    It follows from \eqref{eq:diffDirac} that once $\bar{\theta}$ is an eigenspinor 
    of $D(\bar{g})$, see \eqref{eq:DefEigenOp},
    $\tilde{\Phi}^*\bar{\theta}$ is
    an eigenspinor of $D(r^2\Phi^* g_+)$ with the same eigenvalue. 
    Thus we can calculate the residues of meromorphic spinor valued distribution 
    \begin{align}\label{spindistrcont}
      Res_{\mu=-\frac n2-\lambda-N-1}&\big(M_{\tilde{\Phi}^*\bar{\theta}}(\mu;r)(\phi)\big)\notag\\
      &=Res_{\mu=-\frac n2-\lambda-N-1}\int_X r^\mu<\tilde{\Phi}^*(\bar{\theta}),\phi>Vol(r^2\Phi^* g_+)
    \end{align}
    in two independent ways. On the one hand, integrating $r$ over $[0,\varepsilon)$, using partial integrations
    and extracting residues, (\ref{spindistrcont}) equals to 
    \begin{align*}  
      \int_M e^{(-\lambda-\frac n2)\sigma}&<\theta_0,\slashed{\delta}^+_N(\widehat{h};\lambda)\phi>Vol(\hat{h})\\
      =&\int_M e^{-(\lambda+\frac n2-n)\sigma}<\theta_0,\slashed{\delta}^+_N(\widehat{h};\lambda)\phi>Vol({h}),
    \end{align*}
    because $\frac{Vol(\widehat{h})}{Vol({h})}=e^{n\sigma}$ for the 
    conformal transformation $h\mapsto \widehat{h}=e^{2\sigma}h$.
    On the other hand, equation \eqref{pemetrictrans} and compactness of the support of $\phi$ implies
    \begin{align*}
     &Res_{\mu=-\frac n2-\lambda-N-1}M_{\tilde{\Phi}^*\bar{\theta}}(\lambda;r)(\phi)\\
     =&Res_{\mu=-\frac n2-\lambda-N-1}\int_X r^{\mu +n+1}
       <\Phi^*(\bar{\theta}),\phi>(\Phi^*(r))^{-n-1}Vol(\Phi^*(dr^2+h_r))\\
     =&Res_{\mu=-\frac n2-\lambda-N-1}\int_X \left((\Phi^{-1})^*(r)\right)^{\mu +n+1}
       <\bar{\theta},(\tilde{\Phi}^{-1})^*(\phi)>r^{-n-1}Vol(dr^2+h_r)\\
     =&Res_{\mu=-\frac n2-\lambda-N-1}\int_X \left(\frac{(\Phi^{-1})^*(r)}{r}\right)^{\mu +n+1}r^\mu 
       <\bar{\theta},(\tilde{\Phi}^{-1})^*(\phi)>Vol(dr^2+h_r)\\
     =&\int_M <\theta_0,\slashed{\delta}^+_N(h,\lambda)
        \left(\frac{(\Phi^{-1})^*(r)}{r}\right)^{-\lambda+\frac n2-N}(\tilde{\Phi}^{-1})^*(\phi)>Vol(h).
    \end{align*}
    Since the boundary value $\theta_0$ of $\theta$ was chosen arbitrarily, we get 
    \begin{align*}
       \slashed{\delta}^+_N(\widehat{h},\lambda)=e^{(\frac n2+\lambda-n)\sigma}\circ\slashed{\delta}^+_N(h,\lambda)
          \circ\left(\frac{(\Phi^{-1})^*(r)}{r}\right)^{-\lambda+\frac n2-N}\circ (\tilde{\Phi}^{-1})^*.
    \end{align*}
    For $\lambda\mapsto \lambda+\frac n2-N$ this formula amounts to 
    \begin{align*}
       \slashed{\delta}^+_N(\widehat{h},\lambda+\frac n2-N)
          =e^{(\lambda-N)\sigma}\circ\slashed{\delta}^+_N(h,\lambda+\frac n2-N)
          \circ (\tilde{\Phi}^{-1})^*\circ\left(\frac{\Phi^*(r)}{r}\right)^{\lambda}.
    \end{align*}
    A similar statement holds for $\slashed{\delta}^-_N(h;\lambda+\frac n2-N)$, which completes the proof.
  \end{bew}

  Inspired by \cite[Theorem $6.2.1$]{Juhl}, we construct the first order residue family 
	operators on spinors for general hypersurfaces in conformal manifolds. Let 
  $(X,g)$ be a $(n+1)$-dimensional semi-Riemannian $Spin$-manifold. 
  Consider a hypersurface $\iota:M\to X$ equipped with induced
  metric $h:=\iota^*g$, induced $Spin$-structure 
  and the unit normal vector field $N(g)$. Let $H(g)$ (note the different sign convention 
  for second fundamental form in \cite{Juhl}) be the corresponding 
  mean curvature. Then we define, for $\lambda\in\C$, a family of first order operators 
  \begin{align}
    D_1(X,M;g,\lambda):\Gamma\big(S(X,g)\big)\to&\Gamma\big(S(X,g)|_M\big)\notag\\
    \theta\mapsto& N(g)\cdot\iota^*\slashed{D}^g\theta+(2\lambda+n)\iota^*\nabla_{N(g)}^{g,S}\theta\notag\\
    &+(2\lambda+n)(\lambda-\frac 12)\iota^* H(g)\theta\, .\label{eq:FirstOrderFamilyGeneralHyper}
  \end{align}
  Note that $D_1(X,M;g,\lambda)$ is a polynomial in $\lambda$ of degree 
  two, in contrast to the first order residue family operators. Nevertheless, 
  it satisfies the following conformal transformation law:
  \begin{sat}\label{GeneralHypersurfaceCovariance}
    Let $\widehat{g}=e^{2\sigma}g$ be conformally equivalent to $g$, for $\sigma\in\mathcal{C}^\infty(X)$. Then 
    it holds 
    \begin{align*}
      e^{-(\lambda-1)\iota^*\sigma}D_1(X,M;\widehat{g},\lambda)(e^{\lambda \sigma}\widehat{\theta})
       =\widehat{D_1(X,M;g,\lambda)\theta},
    \end{align*}
    for $\theta\in\Gamma\big(S(X,g)\big)$. Note that $\widehat{\cdot}$ 
    denotes evaluation with respect to $\widehat{g}$. 
  \end{sat}
  \begin{bew}
    The proof based on direct computations. Using 
    \begin{align*}
      \widehat{Y}=&\, e^{-\sigma}Y, \quad \widehat{Y}\hat{\cdot}\widehat{\theta}=\widehat{Y\cdot\theta},
        \quad N(\widehat{g})=e^{-\sigma}N(g),\\
      H(\widehat{g})=&\, e^{-\sigma}\big(H(g)-N(g)(\sigma)\big),\\
      \nabla^{\widehat{g},S}_Y\widehat{\theta}=&\widehat{\nabla_Y^{g,S}\theta}
        -\frac 12\big(\widehat{Y\cdot\grad(\sigma)\cdot\theta}+Y(\sigma)\widehat{\theta}\big),\\
      \slashed{D}^{\widehat{g}}(e^{-\frac{n}{2}\sigma}\widehat{\theta})
        =&\, e^{\frac{n+2}{2}\sigma}\widehat{\slashed{D}^g\theta},
    \end{align*}
    for $Y\in\Gamma(TM)$, $\theta\in\Gamma\big(S(X,g)\big)$, we compute
    \begin{align*}
      N(\widehat{g})\hat{\cdot}\slashed{D}^{\widehat{g}}(e^{\lambda\sigma}\widehat{\theta})
        =&\, e^{(\lambda-1)\iota^*\sigma}\big[\widehat{N(g)\cdot\slashed{D}^g\theta}
        +(\frac n2+\lambda)\widehat{N(g)\cdot\grad(\sigma)\cdot\theta}\big],\\
      (2\lambda+n)\nabla^{\widehat{g},S}_{N(\widehat{g})}(e^{\lambda\sigma}\widehat{\theta})
        =&\, e^{(\lambda-1)\iota^*\sigma}\big[(2\lambda+n)\widehat{\nabla^{g,S}_{N(g)}\theta}
        -(\frac n2+\lambda)\widehat{N(g)\cdot\grad(\sigma)\cdot\theta}\\
      &\, +(2\lambda+n)(\lambda-\frac 12)N(g)(\sigma)\widehat{\theta}\big],\\
     (2\lambda+n)(\lambda-\frac 12)H(\widehat{g})e^{\lambda\sigma}\widehat{\theta}
       =&\, e^{(\lambda-1)\iota^*\sigma}(2\lambda+n)(\lambda-\frac 12)\big[
       H(g)\widehat{\theta}-N(g)(\sigma)\widehat{\theta}\big].
    \end{align*}
    Collecting all terms completes the proof. 
  \end{bew}

  \begin{bem}
     Consider the conformal compactification $(X,\bar{g})$ of the Poincar\'e-Einstein metric $g_+$ associated to $(M,h)$. 
     In this case we have that $\iota^* H=0$ and $N(g)=\partial_r$. Thus, 
     using equation \eqref{eq:HyperDirac} and Example \ref{ExamplesCurved}, we obtain
     \begin{align*}
        D_1(X,M;\bar{g},\lambda)=&\widetilde{\slashed{D}}\iota^*+(2\lambda+n-1)\nabla^{\bar{g},S}_{\partial_r}
          =\slashed{D}^{res}_1(h;\lambda),
     \end{align*}
     hence $D_1(X,M;g,\lambda)$ generalizes the first order residue family operator on spinors. 
  \end{bem}

  The last Proposition indicates the existence of a collection of natural family operators acting on spinors,
  \begin{align*}
     \slashed{D}_N(X, M; g, \lambda): \Gamma\big(S(X,g)\big)\to \Gamma\big(S(X,g)|_{M}\big),
  \end{align*}
  conformally covariant in the sense that 
  \begin{align*}
     e^{-(\lambda-N)(\iota^*\sigma)}\circ \slashed{D}_N(X, M; \widehat{g}, \lambda)\circ e^{\lambda\sigma}=
        \widehat{\slashed{D}_N(X, M; g, \lambda)}
  \end{align*}
  for all $\widehat{g}=e^{2\sigma}g$ and $\lambda\in{\mathbb C}$.
  Notice that in the flat case was the existence of such operators abstractly concluded and their explicit 
  construction was given in \cite{KOSS}, based on the techniques of the classification of homomorphisms of 
  generalized Verma modules.

\section{Factorization identities}\label{FactorizationProp} 

This section presents a complete set of factorization identities for the residue family operators on spinors 
in the flat euclidean case $(\R^{n-1},h)$, while
for a semi-Riemannian $Spin$-manifold $(M^{n},h)$ the factorizations will be 
demonstrated for residue family operators up to order three.

\subsection{Factorization identities - flat case}
  
  Let us start with some low order examples of factorizations:
  \begin{bsp} 
    It directly follows from Example \ref{ExamplesFlat} that low order residue family operators in the flat 
    case obey the following factorizations:
    \begin{align*}
      \slashed{D}^{res}_{1}(h;-\frac{n-2}{2})=&e_n\cdot D_{\mathcal{T}}\iota^*,\quad 
        \slashed{D}^{res}_{1}(h;-\frac{n-1}{2})=e_n\cdot\iota^*(D_{\mathcal{T}}+D_{\mathcal{N}}),\\
      \slashed{D}^{res}_{2}(h;-\frac{n-4}{2})=&e_n\cdot D_{\mathcal{T}}\slashed{D}^{res}_{1}(h;-\frac{n-4}{2})\\
      \slashed{D}^{res}_{2}(h;-\frac{n-1}{2})=
        &-e_n\cdot\slashed{D}^{res}_{1}(h;-\frac{n+1}{2})(D_{\mathcal{T}}+D_{\mathcal{N}}),\\
      \slashed{D}^{res}_{3}(h;-\frac{n-4}{2})=&e_n\cdot D_{\mathcal{T}}^3\iota^*,\quad 
      \slashed{D}^{res}_{3}(h;-\frac{n-3}{2})=e_n\cdot \iota^*(D_{\mathcal{T}}+D_{\mathcal{N}})^3,\\
      \slashed{D}^{res}_{3}(h;-\frac{n-6}{2})=&e_n\cdot D_{\mathcal{T}}\slashed{D}^{res}_{2}(h;-\frac{n-6}{2}),\\
      \slashed{D}^{res}_{3}(h;-\frac{n-1}{2})=
       &e_n\cdot \slashed{D}^{res}_{2}(h;-\frac{n+1}{2})(D_{\mathcal{T}}+D_{\mathcal{N}}).
    \end{align*}
  \end{bsp}

  Now we prove the general form of factorization identities.
  \begin{theorem}\label{FactorizationSmallDirac}
    For $N\geq 1$ and $0\leq M\leq N-1$, the even spinor residue family operators have factorization properties of the form
    \begin{align}
      \slashed{D}^{res}_{2N}(h;2N-\frac n2-M)
        =e_n\cdot D_{\mathcal{T}}^{2M+1} \slashed{D}^{res}_{2N-2M-1}(h;2N-\frac n2-M),\label{eq:EvenFactTang}
    \end{align}
    and for $N\geq 1$ and $0\leq M\leq N$, the odd ones fulfill
    \begin{align}
      \slashed{D}^{res}_{2N+1}(h;2N+1-\frac n2-M)=e_n\cdot D_{\mathcal{T}}^{2M+1}
        \slashed{D}^{res}_{2N-2M}(h;2N+1-\frac n2-M).\label{eq:OddFactTang}
    \end{align}
  \end{theorem}
  \begin{bew}
    In what follows we omit the restriction map $\iota^*$. 
    We start to prove equation \eqref{eq:EvenFactTang}, using explicit formulas for the spinor 
    residue family operators, cf. Theorem \ref{Intertwiner}. First of all, we have
    \begin{align*}
      \slashed{D}^{res}_{2N}&(h;2N-\frac n2-M)
        =(-1)^N\Big[\sum_{j=0}^Na_j^{(N)}(2N-\frac{n+1}{2}-M)D_{\mathcal{T}}^{2j}
        D_{\mathcal{N}}^{2N-2j}\\
      &+2N\sum_{j=0}^{N-1}b_j^{(N-1)}(2N-\frac{n+1}{2}
        -M)D_{\mathcal{T}}^{2j+1}D_{\mathcal{N}}^{2N-2j-1}\Big].
    \end{align*}
    On the other hand, we have
    \begin{align*}
      &\slashed{D}^{res}_{2N-2M-1}(h;2N-\frac n2-M)= \slashed{D}^{res}_{2(N-M-1)+1}(h;2N-\frac n2-M)=\\
      =&(-1)^{N-M}e_n\cdot\Big[\sum_{j=0}^{N-M-1}
         2Nb_j^{(N-M-1)}(2N-\frac{n+1}{2}-M)D_{\mathcal{T}}^{2j}D_{\mathcal{N}}^{2N-2M-2j-1}\\
      &-\sum_{j=0}^{N-M-1}a_j^{(N-M-1)}(2N-\frac{n+1}{2}-M)
         D_{\mathcal{T}}^{2j+1}D_{\mathcal{N}}^{2N-2M-2j-2}\Big].
    \end{align*}
    Multiplying the last formula by $e_n\cdot D_{\mathcal{T}}^{2M+1}=-D_{\mathcal{T}}^{2M+1}e_n\cdot\,$ 
    and shifting the summation index gives
    \begin{align*}
      e_n\cdot D_{\mathcal{T}}^{2M+1}&\slashed{D}^{res}_{2(N-M)-1}(h;2N-\frac n2-M)=\\
      =&(-1)^{N-M}\Big[2N\sum_{j=M}^{N-1}
         b_{j-M}^{(N-M-1)}(2N-\frac{n+1}{2}-M)D_{\mathcal{T}}^{2j+1}D_{\mathcal{N}}^{2N-2j-1}\\
      &-\sum_{j=M+1}^{N}a_{j-M-1}^{(N-M-1)}(2N-\frac{n+1}{2}-M)
         D_{\mathcal{T}}^{2j}D_{\mathcal{N}}^{2N-2j}\Big].
    \end{align*}
    This equals to $\slashed{D}^{res}_{2N}(h;2N-\frac n2-M)$ provided we have for all $j=M+1,\ldots, N$
    \begin{align}
       (-1)^Na_j^{(N)}(2N-\frac{n+1}{2}-M)
         =(-1)^{N-M-1}a_{j-M-1}^{(N-M-1)}(2N-\frac{n+1}{2}-M),\label{eq:condition1}
    \end{align}
    and for all $j=0,\ldots,M$ we have $a_j^{(N)}(2N-\frac{n+1}{2}-M)=0$, as well as for 
    all $j=M,\ldots, N-1$ we have 
    \begin{align}
       (-1)^Nb_j^{(N-1)}(2N-\frac n2-M-\frac 12)
         =(-1)^{N-M}b_{j-M}^{(N-M-1)}(2N-\frac n2-M-\frac 12),\label{eq:condition2}
    \end{align}
    and for all $j=0,\ldots,M-1$ holds $b_j^{(N-1)}(2N-\frac{n+1}{2}-M)=0$. Using definition \eqref{eq:DefA}, 
    it directly follows that equation \eqref{eq:condition1} is equivalent to
    \begin{align*}
      \frac{N!}{j!}&\prod_{k=j}^{N-1}(2k-2M)=\frac{(N-M-1)!}{(j-M-1)!}\prod_{k=j-M-1}^{N-M-2}(2k+2M+4).
    \end{align*}
    This is easily verified, and we also get from definition \eqref{eq:DefA} the
    triviality of coefficients $a_j^{(N)}(2N-\frac{n+1}{2}-M)=0$ for 
    all $j=0,\ldots,M$. One can also prove, using definition \eqref{eq:DefB}, similar statement for $b_j^{(N-1)}$.

    Now we discuss equation \eqref{eq:OddFactTang}. From Theorem \ref{Intertwiner} we have
    \begin{align*}
      \slashed{D}^{res}_{2N+1}&(h;2N+1-\frac n2-M)=\\
      =&(-1)^{N+1}e_n\cdot\Big[2(N-M)
        \sum_{j=0}^{N} b_j^{(N)}(2N+\frac{1-n}{2}-M)D_{\mathcal{T}}^{2j}D_{\mathcal{N}}^{2N-2j+1}\\
      &-\sum_{j=0}^N a_j^{(N)}(2N+\frac{1-n}{2}-M)D_{\mathcal{T}}^{2j+1}D_{\mathcal{N}}^{2N-2j}\Big],
    \end{align*}
    as well as
    \begin{align*}
      &\slashed{D}^{res}_{2(N-M)}(h;2N+1-\frac n2-M)=\\
      =&(-1)^{N-M}\Big[\sum_{j=0}^{N-M}a_{j}^{(N-M)}(2N+\frac{1-n}{2}-M)
       D_{\mathcal{T}}^{2j}D_{\mathcal{N}}^{2(N-M)-2j}\\
      &+2(N-M)\sum_{j=0}^{N-M-1}b_j^{(N-M-1)}(2N+\frac{1-n}{2}-M)
       D_{\mathcal{T}}^{2j+1}D_{\mathcal{N}}^{2(N-M)-2j-1}\Big].      
    \end{align*}
    We multiply the last formula by $e_n\cdot D_T^{2M+1}$ and shift the summation index, resulting in 
    \begin{align*}
      e_n\cdot & D_T^{2M+1}\slashed{D}^{res}_{2(N-M)}(h;2N+1-\frac n2-M)=\\
      =&(-1)^{N-M}e_n\cdot\Big[\sum_{j=M}^{N}a_{j-M}^{(N-M)}(2N+\frac{1-n}{2}-M)
       D_{\mathcal{T}}^{2j+1}D_{\mathcal{N}}^{2N-2j}\\
      &+2(N-M)\sum_{j=M+1}^{N}b_{j-M-1}^{(N-M-1)}(2N+\frac{1-n}{2}-M)
       D_{\mathcal{T}}^{2j}D_{\mathcal{N}}^{2N-2j+1}\Big].
    \end{align*}
    This equals to $\slashed{D}^{res}_{2N+1}(h;2N+1-\frac n2-M)$ provided we have for all $j=M,\ldots, N$
    \begin{align*}
      (-1)^Na_j^{(N)}(2N+\frac{1-n}{2}-M)=(-1)^{N-M}a^{(N-M)}_{j-M}(2N+\frac{1-n}{2}-M),
    \end{align*}
    and for all $j=0,\ldots,M-1$ we have $a_j^{(N)}(2N+\frac{1-n}{2}-M)=0$, as well 
    as for all $j=M+1,\ldots, N$ we have
    \begin{align*}
      (-1)^{N}b_j^{(N)}(2N+\frac{1-n}{2}-M)=(-1)^{N-M-1}b^{(N-M-1)}_{j-M-1}(2N+\frac{1-n}{2}-M),
    \end{align*}
    and for all $j=0,\ldots,M$ we have $b_j^{(N)}(2N+\frac{1-n}{2}-M)=0$. 
    These properties are easily checked from definitions 
    \eqref{eq:DefA} and \eqref{eq:DefB}, and this completes the proof.
  \end{bew}

  The next theorem states that the residue family operators afford factorization properties from the right 
  by powers of the Dirac operator $D_{\mathcal{T}}+D_{\mathcal{N}}$ on $\R^n$. 
  Since this proof is more difficult than the left factorization, 
  we start with two preparatory lemmas.
  \begin{lem}
    For $M\in\N_0$, we have 
    \begin{align*}
      (D_{\mathcal{T}}+D_{\mathcal{N}})^{2M+1}
        =&\sum_{l=0}^{2M+1}\begin{pmatrix} M\\ \lfloor \frac l2\rfloor\end{pmatrix} 
           D_{\mathcal{T}}^lD_{\mathcal{N}}^{2M+1-l}\\
      =&\sum_{l=0}^M\begin{pmatrix} M\\ l\end{pmatrix} D_{\mathcal{T}}^{2l}D_{\mathcal{N}}^{2M+1-2l}
       +\sum_{l=0}^M\begin{pmatrix} M\\ l\end{pmatrix} D_{\mathcal{T}}^{2l+1}D_{\mathcal{N}}^{2M-2l},
    \end{align*}
    where $\lfloor k \rfloor$ denotes the integer part of $k$.
  \end{lem}
  \begin{bew}
    The proof goes by induction. The case $M=0$ holds for trivial reasons, and we assume the statement holds for $M$. 
    Then we have for $M+1$:
    \begin{align*}
      (D_{\mathcal{T}}+D_{\mathcal{N}})^{2(M+1)+1}=&(D_{\mathcal{T}}+D_{\mathcal{N}})^{2M+1}
           (D_{\mathcal{T}}+D_{\mathcal{N}})^2\\
      =&\sum_{l=0}^{2M+1}\begin{pmatrix} M\\ \lfloor \frac l2\rfloor\end{pmatrix} 
        D_{\mathcal{T}}^lD_{\mathcal{N}}^{2M+1-l}(D_{\mathcal{T}}+D_{\mathcal{N}})^2\\
	      =&\sum_{l=2}^{2M+3}\begin{pmatrix} M\\ \lfloor \frac{l-2}{2}\rfloor\end{pmatrix} 
        D_{\mathcal{T}}^lD_{\mathcal{N}}^{2M+3-l}
       +\sum_{l=0}^{2M+1}\begin{pmatrix} M\\ \lfloor \frac l2\rfloor\end{pmatrix} 
        D_{\mathcal{T}}^lD_{\mathcal{N}}^{2M+3-l}\\
      =&\sum_{l=0}^{2(M+1)+1}\begin{pmatrix} M+1\\ \lfloor \frac l2\rfloor\end{pmatrix}D_{\mathcal{T}}^l
          D_{\mathcal{N}}^{2(M+1)+1-l},
    \end{align*} 
    which completes the induction step. Decomposing the sum into even and odd values of 
    the summation index $l$ gives the desired last statement. 
  \end{bew}

  The next lemma treats several combinatorial identities involving binomial coefficients and Pochhammer symbols. 
  We recall the definition of Pochhammer symbol $(a)_l=a(a+1)\ldots(a+l-1)$ and $(a)_0:=1$, for $a\in\R$ 
  and $l\in\N_0$. 
  \begin{lem}\label{PochhammerIdentities}
    \begin{enumerate}
     \item
      For all $j,M,N\in{\mathbb N}$ such that $1\leq M\leq N-M-1$ and $M\leq j\leq N-M-1$,
      the identity
      \begin{align} 
        (N-M)&_M(N-M+\frac{1}{2})_{M}=\notag\\
        =& \sum_{l=0}^M(-1)^{l}{M\choose l}(j-l+1)_l(N-M-j+l)_{M-l}\times\notag\\
        &\times (N-M-j+l+\frac{1}{2})_{M-l}(j-l-2N+M+\frac{3}{2})_{l},\label{eq:binform}
      \end{align}
      holds true.
     \item
      For all $j,M,N\in{\mathbb N}$ such that $1\leq M\leq N-M-1$ and $M\leq j\leq N-M-1$,
      the identity
      \begin{align*} 
        (N-M+\frac 12)&_M(N-M+1)_{M}=\\
        &\sum_{l=0}^M(-1)^{l}{M\choose l}(j-l+1)_l(N-M-j+l+\frac 12)_{M-l}\times\\
        &\times (N-M-j+l+1)_{M-l}(j-l-2N+M+\frac{1}{2})_{l}
      \end{align*}
      holds true. 
     \item
      For all $j,M,N\in{\mathbb N}$ such that $1\leq M\leq N-M$ and $M+1\leq j\leq N-M$, the identity
      \begin{align*} 
        (N-M+1)&_M(N-M+\frac{3}{2})_{M}=\\
        =&\sum_{l=0}^M(-1)^{l}{M\choose l}(j-l+1)_l(N-M-j+l+1)_{M-l}\times\\
        &\times (N-M-j+l+\frac{3}{2})_{M-l}(j-l-2N+M-\frac{1}{2})_{l},
      \end{align*}
      holds true.
    \end{enumerate}
  \end{lem}

  \begin{bew}
    We shall prove only the first statement, since the second and third are proved analogously. 
    Recall the definition of a generalized hypergeometric function 
    ${}_pF_{q}(a_1,a_2,\dots ,a_p;b_1,b_2,\dots ,b_q;z)$ defined by 
    \begin{align}
     {}_pF_{q}=
      \begin{pmatrix}
       a_1, & a_2, & \dots  &,\, a_p & \\
        &    &       &    &; z \\
       b_1, & b_2, & \dots  &,\, b_q &
      \end{pmatrix}
      =\sum_{n=0}^\infty \frac{(a_1)_n(a_2)_n\dots (a_p)_n}{(b_1)_n(b_2)_n\dots (b_q)_n}\frac{z^n}{n!}.
    \end{align} 
    A direct comparison converts our binomial expression into balanced ${}_3F_{2}$ series, 
    \begin{align*}
      &{}_3F_{2}
        \begin{pmatrix}
          -j, &2N-M-j+\frac{1}{2},& -M \\
          &    &       &    ; 1 \\
         N-M-j, & N-M-j+\frac{1}{2} &  & \\
        \end{pmatrix} \times \\
      &\times 4^{-M}(-2j-2M+2N)_{2M} =\sum_{l=0}^\infty(-1)^{l}{M\choose l}(j-l+1)_l\times\\
      &4^{l-M}\times (2N-2M-2j+2l)_{2(M-l)}(j-l-2N+M+\frac{3}{2})_{l}, 
    \end{align*}
    which terminates the summation index by $l=M$. Now we use the Pfaff-Saalschutz summation formula,
    cf. \cite[Appendix III.2]{sl},
    \begin{align}
      {}_3F_{2}
      \begin{pmatrix}
        a, & b, & -n  \\
         &    &       &    ; 1 \\
        c, & 1+a+b-c-n &  & \\
      \end{pmatrix} 
      =\frac{(-a+c)_n(-b+c)_n}{(c)_n(-a-b+c)_n},
    \end{align}
    for $n\in{\mathbb N}$. This converts the hypergeometric terminating expansion 
    for ${}_3F_{2}$ into required binomial identity and the proof is complete.
  \end{bew}

  \begin{bem}
    The proof of Lemma \ref{PochhammerIdentities} 
    was communicated to the authors by Ch. Krattenthaler. He also informed us 
    that there is a Gosper-Zeilberger algorithm, cf. 
    \cite{ps} and references therein for its implementation, allowing to find a recurrence 
    for binomial sums. In particular, denoting the right hand side of equation \eqref{eq:binform} by 
    $S[N]$, one finds the functional recurrence relation
    \begin{align}
      N(1+2N)S[N]-(-1+2M-2N)(M-N)S[N+1]=0.
    \end{align}
  \end{bem}
  Now we can prove the second set of factorizations.
  \begin{theorem}\label{FactorziationBicDirac}
    For $N\geq 1$ and $0\leq M\leq N-1$, the even spinor residue family operators have the right factorization property
    \begin{align}
      \slashed{D}^{res}_{2N}(h;\frac{1-n}{2}+M)=-e_n\cdot \slashed{D}^{res}_{2(N-M-1)+1}(h;-\frac{1+n}{2}-M)
        (D_{\mathcal{T}}+D_{\mathcal{N}})^{2M+1},\label{eq:EvenFactBig}
    \end{align}
    while for $N\geq 1$ and $0\leq M\leq N$, the odd ones fulfill
    \begin{align}
      \slashed{D}^{res}_{2N+1}(h;\frac{1-n}{2}+M)=e_n\cdot 
        \slashed{D}^{res}_{2N-2M}(h;-\frac{n+1}{2}-M)
          (D_{\mathcal{T}}+D_{\mathcal{N}})^{2M+1}.\label{eq:OddFactBig}
    \end{align}
  \end{theorem}
  \begin{bew}
    It follows from Theorem \ref{Intertwiner} that
    \begin{align*}
      \slashed{D}^{res}_{2N}(h;\frac{1-n}{2}+M)
       =&(-1)^N\Big[\sum_{j=0}^Na_j^{(N)}(-\frac n2+M)D_{\mathcal{T}}^{2j}D_{\mathcal{N}}^{2N-2j}\\
      &+2N\sum_{j=0}^{N-1}b_j^{(N-1)}(-\frac n2+M)D_{\mathcal{T}}^{2j+1}D_{\mathcal{N}}^{2N-2j-1}\Big],
    \end{align*}
    and
    \begin{align*}
      &\slashed{D}^{res}_{2(N-M-1)+1}(h;-\frac{1-n}{2}-M)=\\
      =&(-1)^{N-M}e_n\cdot\Big[-\sum_{j=0}^{N-M-1}a_j^{(N-M-1)}(-\frac n2-M-1)
           D_{\mathcal{T}}^{2j+1}D_{\mathcal{N}}^{2(N-M-1)-2j}\\
      &-(2N+1)\sum_{j=0}^{N-M-1}b_j^{(N-M-1)}(-\frac n2-M-1)
        D_{\mathcal{T}}^{2j}D_{\mathcal{N}}^{2(N-M-1)-2j+1}\Big].
    \end{align*}
    Multiplying the last formula from the left by $-e_n\cdot$ and from the right by 
    $(D_{\mathcal{T}}+D_{\mathcal{N}})^{2M+1}$ gives
    \begin{align}
      &-e_n\cdot\slashed{D}^{res}_{2(N-M-1)+1}(h;-\frac{1-n}{2}-M)
       (D_{\mathcal{T}}+D_{\mathcal{N}})^{2M+1}=\notag\\
      =&(-1)^{N-M}\Big[-(2N+1)\sum_{j=0}^{N-M-1}b_j^{(N-M-1)}(-\frac n2-M-1)
         D_{\mathcal{T}}^{2j}D_{\mathcal{N}}^{2(N-M-1)-2j+1}\notag\\
      &-\sum_{j=0}^{N-M-1}a_j^{(N-M-1)}(-\frac n2-M-1)
         D_{\mathcal{T}}^{2j+1}D_{\mathcal{N}}^{2(N-M-1)-2j}\Big]
        (D_{\mathcal{T}}+D_{\mathcal{N}})^{2M+1}\notag\\
      =&(-1)^{N-M}\Big[
       -(2N+1)\sum_{j=0}^{N-M-1}\sum_{l=0}^Mb_j^{(N-M-1)}(-\frac n2-M-1)
        \begin{pmatrix} M\\ l\end{pmatrix}D_{\mathcal{T}}^{2j+2l}D_{\mathcal{N}}^{2N-2j-2l}\notag\\
      &-\sum_{j=0}^{N-M-1}\sum_{l=0}^M a_j^{(N-M-1)}(-\frac n2-M-1)
       \begin{pmatrix} M\\ l\end{pmatrix}D_{\mathcal{T}}^{2j+2l+2}D_{\mathcal{N}}^{2N-2j-2l-2}\notag\\
      &+(2N+1)\sum_{j=0}^{N-M-1}\sum_{l=0}^M b_j^{(N-M-1)}(-\frac n2-M-1)
        \begin{pmatrix} M\\ l\end{pmatrix}D_{\mathcal{T}}^{2j+2l+1}D_{\mathcal{N}}^{2N-2j-2l-1}\notag\\
      &-\sum_{j=0}^{N-M-1}\sum_{l=0}^M a_j^{(N-M-1)}(-\frac n2-M-1)
        \begin{pmatrix} M\\ l\end{pmatrix}D_{\mathcal{T}}^{2j+2l+1}D_{\mathcal{N}}^{2N-2j-2l-1}\Big]\notag\\
      =&(-1)^{N-M}\Bigg[
       -(2N+1)\sum_{l=0}^Mb_0^{(N-M-1)}(-\frac n2-M-1)
       \begin{pmatrix} M\\ l\end{pmatrix}D_{\mathcal{T}}^{2l}D_{\mathcal{N}}^{2N-2l}\notag\\
      &-\sum_{l=N-M}^N a^{(N-M-1)}_{N-M-1}(-\frac n2-M-1)
       \begin{pmatrix} M\\ l-N+M\end{pmatrix}D_{\mathcal{T}}^{2l}D_{\mathcal{N}}^{2N-2l}\notag\\
      &-\sum_{j=1}^{N-M-1}\sum_{l=0}^M\Big[(2N+1)b_j^{(N-M-1)}(-\frac n2-M-1)\notag\\
      &\quad\quad\quad\quad\quad\quad\quad+a_{j-1}^{(N-M-1)}(-\frac n2-M-1)\Big]
       \begin{pmatrix} M \\ l\end{pmatrix}D_{\mathcal{T}}^{2j+2l}D_{\mathcal{N}}^{2N-2j-2l}\label{eq:sum1}\\
      &+\sum_{j=0}^{N-M-1}\sum_{l=0}^M\Big[ (2N+1)b_j^{(N-M-1)}(-\frac n2-M-1)\notag \nonumber \\ 
      &\quad\quad\quad\quad\quad\quad\quad-a_j^{(N-M-1)}(-\frac n2-M-1)\Big]
        \begin{pmatrix} M \\ l\end{pmatrix}
        D_{\mathcal{T}}^{2j+2l+1}D_{\mathcal{N}}^{2N-2j-2l-1}\Bigg].\label{eq:sum2}
    \end{align}
    Let us further assume $M\leq N-M-1$, the case $N-M-1< M$ is analogous and follows from the 
    first by transformation $M\to N-M-1$. We introduce the abbreviations 
    \begin{align}
      c_j:=&(2N+1)b_j^{(N-M-1)}(-\frac n2-M-1)+a_{j-1}^{(N-M-1)}(-\frac n2-M-1),\\
      d_j:=&(2N+1)b_j^{(N-M-1)}(-\frac n2-M-1)-a_{j}^{(N-M-1)}(-\frac n2-M-1),
    \end{align}
    for which two sums in equation \eqref{eq:sum1} and \eqref{eq:sum2} can be reformulated as
    \begin{align*}
      &\sum_{j=1}^{N-M-1}\sum_{l=0}^Mc_j\begin{pmatrix} M\nonumber \\ l\end{pmatrix}
        D_{\mathcal{T}}^{2j+2l}D_{\mathcal{N}}^{2N-2j-2l}=\nonumber \\
      =&\sum_{i=1}^M\sum_{l=0}^{i-1}c_{i-l}\begin{pmatrix} M\nonumber \\ l\end{pmatrix}
        D_{\mathcal{T}}^{2i}D_{\mathcal{N}}^{2N-2i}
        +\sum_{i=M+1}^{N-M-1}\sum_{l=0}^{M}c_{i-l}\begin{pmatrix} M\nonumber \\ l\end{pmatrix}
        D_{\mathcal{T}}^{2i}D_{\mathcal{N}}^{2N-2i}\\
      &+\sum_{i=N-M}^{N-1}\sum_{l=i-N+M+1}^{M}c_{i-l}\begin{pmatrix} M\\ l\end{pmatrix}
        D_{\mathcal{T}}^{2i}D_{\mathcal{N}}^{2N-2i},
    \end{align*}
    respectively,
    \begin{align*}
      &\sum_{j=0}^{N-M-1}\sum_{l=0}^Md_j\begin{pmatrix} M\nonumber \\ l\end{pmatrix}
        D_{\mathcal{T}}^{2j+2l+1}D_{\mathcal{N}}^{2N-2j-2l-1}=\nonumber \\
      =&\sum_{i=0}^{M-1}\sum_{l=0}^{i}d_{i-l}\begin{pmatrix} M\nonumber \\ l\end{pmatrix}
       D_{\mathcal{T}}^{2i+1}D_{\mathcal{N}}^{2N-2i-1}
        +\sum_{i=M}^{N-M-1}\sum_{l=0}^{M}d_{i-l}\begin{pmatrix} M\nonumber \\ l\end{pmatrix}
       D_{\mathcal{T}}^{2i+1}D_{\mathcal{N}}^{2N-2i-1}\nonumber\\
      &+\sum_{i=N-M}^{N-1}\sum_{l=i-N+M+1}^{M}d_{i-l}
       \begin{pmatrix} M\\ l\end{pmatrix}D_{\mathcal{T}}^{2i+1}D_{\mathcal{N}}^{2N-2i-1}.
    \end{align*}
    Then equation \eqref{eq:EvenFactBig} is equivalent, when comparing the coefficients by
    $D_{\mathcal{T}}^kD_{\mathcal{N}}^{2N-k}$, to the following system of equations: 
    for even $k$ we get
    \begin{align}
      (-1)^Na_0^{(N)}(-\frac n2+M)=&(-1)^{N-M-1}(2N+1)b_0^{(N-M-1)}(-\frac n2-M-1),\label{eq:check1}\\
      (-1)^Na_j^{(N)}(-\frac n2+M)=&(-1)^{N-M-1}\sum_{l=0}^j
        \begin{pmatrix} M\\ l\end{pmatrix}c_{j-l},\quad\forall\; 1\leq j\leq M,\label{eq:check2}\\
      (-1)^Na_j^{(N)}(-\frac n2+M)=&(-1)^{N-M-1}\sum_{l=0}^M
        \begin{pmatrix} M\\ l\end{pmatrix}c_{j-l},\quad\forall\; M+1\leq j\leq N-M-1,\label{eq:check3}\\
      (-1)^Na_j^{(N)}(-\frac n2+M)=&(-1)^{N-M-1}\sum_{l=j-N+M}^M
        \begin{pmatrix} M\\ l\end{pmatrix}c_{j-l},\quad\forall\; N-M\leq j\leq N-1,\label{eq:check4}
    \end{align}
    whereas for odd $k$ we have
    \begin{align}
      (-1)^{N}2Nb_j^{(N-1)}(-\frac n2+M)=&(-1)^{N-M}\sum_{l=0}^j\begin{pmatrix} M\\ l\end{pmatrix}d_{j-l},
        \quad \forall\; 0\leq j\leq M-1,\label{eq:check5}\\
      (-1)^{N}2Nb_j^{(N-1)}(-\frac n2+M)=&(-1)^{N-M}\sum_{l=0}^M\begin{pmatrix} M\\ l\end{pmatrix}d_{j-l},
        \quad \forall\; M\leq j\leq N-M-1,\label{eq:check6}\\
      (-1)^{N}2Nb_j^{(N-1)}(-\frac n2+M)=&(-1)^{N-M}\sum_{l=j-N+M+1}^M
        \begin{pmatrix} M\\ l\end{pmatrix}d_{j-l},
        \quad \forall\; N-M\leq j\leq N-1.\label{eq:check7}
    \end{align}
    Based on  definitions \eqref{eq:DefA} and \eqref{eq:DefB}, an easy computation proves equation \eqref{eq:check1}. 
    Equations \eqref{eq:check2} up to \eqref{eq:check7} can be check in a quite uniform way. First of all, 
    definitions \eqref{eq:DefA} and \eqref{eq:DefB} allow us to compute
    \begin{align*}
      (-&1)^{N-M-1}c_{j-l}=(2N+1)b_{j-l}^{(N-M-1)}(-\frac n2-M-1)+a_{j-l-1}^{(N-M-1)}(-\frac n2-M-1)\\
      =&(2N+1)\frac{(N-M-1)!(-4)^{N-M-1-j+l}}{(j-l)!(2N-2M-2j+2l-1)!}(j-l-2N+M+\frac 12)_{N-M-1-j+l}\\
      &+\frac{(N-M-1)!(-4)^{N-M-j+l}}{(j-l-1)!(2N-2M-2j+2l)!}(j-l-2N+M+\frac 12)_{N-M-j+l}\\
      =&\frac{(N-M-1)!(-4)^{N-M-1-j+l}}{(j-l-1)!(2N-2M-2j+2l-1)!}
        \Big[ \frac{(-2)}{(j-l)}+\frac{(-4)}{(2N-2M-2j+2l)}\Big]\times\\
      &\times(j-l-2N+M+\frac 12)_{N-M-j+l}\\
     =&\frac{(N-M-1)!(-4)^{N-M-1-j+l}}{(j-l-1)!(2N-2M-2j+2l-1)!}
        \Big[ \frac{(-4)(N-M)}{(j-l)(2N-2M-2j+2l)}\Big]\times\\
      &\times(j-l-2N+M+\frac 12)_{N-M-j+l}\\
      =&\frac{(N-M)!(-4)^{N-M-j+l}}{(j-l)!(2N-2M-2j+2l)!}(j-l-2N+M+\frac 12)_{N-M-j+l},
    \end{align*}
    and similarly one gets
    \begin{align*}
      (-&1)^{N-M-1}d_{j-l}=(2N+1)b_j^{(N-M-1)}(-\frac n2-M-1)-a_{j}^{(N-M-1)}(-\frac n2-M-1)\\
      =&2N\frac{(N-M-1)!(-4)^{N-M-j+l-1}}{(j-l)!(2N-2M-2j+2l-1)!}(j-l-2N+M+\frac 32)_{N-M-j+l-1}.
    \end{align*}
    Thus, for example, equation \eqref{eq:check3} is equivalent to
    \begin{align*}
      &\frac{N!(-4)^{N-j}}{j!(2N-2j)!}(j-2N+M+\frac 12)_{N-j}=\\
      =&\sum_{l=0}^M {M\choose l}\frac{(N-M)!(-4)^{N-M-j+l}}{(j-l)!(2N-2M-2j+2l)!}(j-l-2N+M+\frac 12)_{N-M-j+l},
    \end{align*}
    which can be reduced, using Pochhammer identities, to  
    \begin{align*}
      (N-M&+1)_M(N-M+\frac 12)_M=\sum_{l=0}^M(-1)^l(4)^{l-M}{M\choose l}\times\\
      &\times (j-l-2N+M+\frac 12)_{l}(2N-2M-2j+2l+1)_{2(M-l)}(j-l+1)_l\\
      =&\sum_{l=0}^M(-1)^l{M\choose l} (j-l-2N+M+\frac 12)_{l}\times\\
      &\times (N-M-j+l+\frac 12)_{M-l}(N-M-j+l+1)_{M-l}(j-l+1)_l.
    \end{align*}
    Lemma \ref{PochhammerIdentities} allows us to conclude that the last equality is true, and 
    therefore equation \eqref{eq:check3} holds. Similarly, 
    one observes that equations \eqref{eq:check2} and \eqref{eq:check4} follow from Lemma \ref{PochhammerIdentities}, 
    due to the range of the index $j$ making trivial some of the Pochhammer symbols and allowing to rearrange 
    the summation array accordingly. 
    Let us proceed to the proof of equation \eqref{eq:check6}. This amounts to check
    \begin{align}
      2N\frac{(N-1)!(-4)^{N-j-1}}{j!(2N-2j-1)!}&(j-2N+M+\frac 32)_{N-j}=\nonumber\\
      =&2N\sum_{l=0}^M {M\choose l}\frac{(N-M-1)!(-4)^{N-M-j+l-1}}{(j-l)!(2N-2M-2j+2l-1)!}\times\nonumber\\
      &\times   (j-l-2N+M+\frac 32)_{N-M-j+l-1},
    \end{align}
    which is equivalent, using Pochhammer identities, to
    \begin{align}
       (N-M&)_M(N-M+\frac 12)_M=\sum_{l=0}^M(-1)^l{M\choose l}(j-2N+M+\frac 32)_l\times\nonumber\\
       &\times (N-M-j+l)_{M-l}(N-M-j+l+\frac 12)_{M-l}(j-l+1)_l,
    \end{align}
    and this holds true due to Lemma \ref{PochhammerIdentities}. Again, 
    equations \eqref{eq:check5} and \eqref{eq:check7} can 
    be checked analogously using Lemma \ref{PochhammerIdentities} and taking into account the range of the index $j$. 

    Now we pass to the second set of factorizations. Based on abbreviations 
    \begin{align}
      c_j:=&-2(N-M)b_{j-1}^{(N-M-1)}(-\frac n2-M-1)+a_{j}^{(N-M)}(-\frac n2-M-1),\\
      d_j:=&2(N-M)b_j^{(N-M-1)}(-\frac n2-M-1)-a_{j}^{(N-M)}(-\frac n2-M-1),
    \end{align}
    equation \eqref{eq:OddFactBig} is equivalent to the system of equalities
    \begin{align}
      (-1)^{N}(2N-2M-1)b_0^{(N)}&(-\frac n2+M)=(-1)^{N-M}a_0^{(N-M)}(-\frac n2-M-1),\label{eq:check7a}\\
      (-1)^{N}(2N-2M-1)b_j^{(N)}&(-\frac n2+M)\notag\\
      =&\, (-1)^{N-M}\sum_{l=0}^j{M\choose l} c_{j-l}, \forall\; 1\leq j\leq M,\label{eq:check8}\\
      (-1)^{N}(2N-2M-1)b_0^{(N)}&(-\frac n2+M)\notag\\
      =&\, (-1)^{N-M}\sum_{l=0}^j{M\choose l} c_{j-l},\forall\; M+1\leq j\leq N-M,\label{eq:check9}\\
      (-1)^{N}(2N-2M-1)b_0^{(N)}&(-\frac n2+M)\notag\\
      =&\, (-1)^{N-M}\sum_{l=0}^j{M\choose l} c_{j-l},\forall\; N-M+1\leq j\leq N,\label{eq:check10}
    \end{align}  
    and
    \begin{align}
      (-1)^{N}a_j^{(N)}&(-\frac n2+M)=(-1)^{N-M}\sum_{l=0}^j {M\choose l} d_{j-l},
        \forall\; 0\leq j\leq M-1,\label{eq:check11}\\
      (-1)^{N}a_j^{(N)}&(-\frac n2+M)=(-1)^{N-M}\sum_{l=0}^M {M\choose l} d_{j-l},
        \forall\; M\leq j\leq N-M-1,\label{eq:check12}\\
      (-1)^{N}a_j^{(N)}&(-\frac n2+M)=(-1)^{N-M}\sum_{l=j-N+M}^M {M\choose l} d_{j-l},
        \forall\; N-M\leq j\leq N-1,\label{eq:check13}
    \end{align}
    By similar arguments as for the factorizations of $\slashed{D}^{res}_{2N}(h;\lambda)$ it follows that 
    equations \eqref{eq:check7a} - \eqref{eq:check13} holds. 
    This completes the proof of the theorem. 
  \end{bew}

  A direct consequence of the last theorem (the case $M=N$) together with the conformal 
  covariance of residue family operators, cf. Theorem \ref{ConformalTrafoLawResidueFamily}, is a 
  well known observation:
  \begin{kor}
    Let $N\in\N_0$. The $(2N+1)$-th power of the Dirac operator on the flat euclidean space is conformally covariant. 
  \end{kor} 
  
  \begin{bem}\label{WhyNoIdentification}
  We finish by several remarks on the target space of residue family operators
    \begin{align*}
      \slashed{D}^{res}_{N}(h;\lambda):\Gamma\big(S(\R^n,\bar{g})\big)\to\Gamma(S(\R^{n},\bar{g})|_{x_n=0}).
    \end{align*}
    In the case when the identification map $\Xi$, see Remark 
    \ref{BundleIdentificationCurved}, is included into the definition 
    of residue family operators, the first and second order family operators can be written as
    \begin{align*}
      \slashed{D}^{res}_{1}(h;\lambda)
       =&D_{\mathcal{T}}\circ\Xi\circ\iota^*-(2\lambda+n-2)\Xi\circ\iota^*e_n\cdot D_{\mathcal{N}},\\
      \slashed{D}^{res}_{2}(h;\lambda)
       =&D_{\mathcal{T}}^2\circ\Xi\circ\iota^*-2D_{\mathcal{T}}\circ\Xi\circ\iota^*e_n\cdot D_{\mathcal{N}}
         -(2\lambda+n-4)\Xi\circ\iota^*D_{\mathcal{N}}^2,
    \end{align*}
    hence they are maps from $\Gamma\big(S(\R^n,\bar{g})\big)$ to $\Gamma\big(S(\R^{n-1},h)\big)$.
    The simultaneous change of $\Xi:=\Xi^-$, cf. Remark \ref{BundleIdentificationCurved}, then leads to a partial 
    change of signs. Consequently, one can not prove the factorization properties in general: some 
    factorizations still survive, but for example for $\slashed{D}^{res}_{2}$ we have 
    \begin{align*}
       \slashed{D}^{res}_{2}(h;-\frac{n-1}{2})=
          D_{\mathcal{T}}^2\circ\Xi\circ\iota^*-2D_{\mathcal{T}}\circ\Xi\circ\iota^*e_n\cdot D_{\mathcal{N}}
         +3\Xi\circ\iota^*D_{\mathcal{N}}^2,
    \end{align*}
    whereas 
    \begin{align*}
      \slashed{D}^{res}_{1}(h;-\frac{n+1}{2})(D_{\mathcal{T}}+D_{\mathcal{N}})
        =&D_{\mathcal{T}}^2\circ\Xi\circ\iota^*e_n\cdot-2D_{\mathcal{T}}\circ\Xi\circ\iota^*D_{\mathcal{N}}\\
        &+3\Xi\circ\iota^*e_n\cdot D_{\mathcal{N}}^2.
    \end{align*}
    Observe that there is no way to get an identity between the two former formulas. 
    Factorization identities 
    hold true if one considers $S(\R^n,\bar{g})|_{x_n=0}$ as a target space for residue family operators, i.e., omitting 
    the identification map $\Xi$. 
  \end{bem}

\subsection{Factorization identities of low order - curved case}

  We start with a general description of the geometry of embedded hypersurfaces and related techniques 
  associated to the geometry of $Spin$-structures, cf. \cite{Bures, BGM}, used to materialize the factorization identities 
  of low order residue 
  families operators on spinors. 
 
  Let $\iota:M^n\hookrightarrow Z^{n+1}$ be an embedded hypersurface or boundary of a semi-Riemannian 
  manifold $(Z,g)$ with induced metric $h:=\iota^*(g)$ on $M$. The restriction and projection 
  of the Levi-Civita connection $\nabla^g$ on $Z$ to $TM$ is the Levi-Civita
  connection $\nabla^h$ on $(M,h)$. Assuming the existence of a unit normal vector field $\nu$ to $M$
  in $Z$, the second fundamental form of $\iota:M\hookrightarrow Z$ is a scalar valued symmetric 
  $2$-tensor 
  \begin{align*}
    II:\Gamma(TM)\times \Gamma(TM)&\to \mathcal{C}^{\infty}(M),\\
    X,Y&\mapsto g(\nabla^g_XY,\nu).
  \end{align*}
  The mean curvature of $\iota:M\hookrightarrow Z$ is the normalized $h$-trace of the second fundamental 
  form $II$, i.e., in a local $h$-orthonormal frame $(s_1,\dots ,s_n)$ we have
  \begin{align*}
    H=\frac 1n\sum_{i=1}^{n}\varepsilon_i II(s_i,s_i).
  \end{align*}
  For our purposes it is sufficient to consider the concept of generalized cylinders, i.e.,  
  $(Z^{n+1}=M^n\times I,g=dr^2+h_r)$, where $h_r$ is a $1$-parameter family of metrics on $M$ and $\partial_r$ 
  denotes the coordinate vector field, 
  cf. \cite{BGM}. For every $r\in I$ let $M_r:=(M\times \{r\}, h_r)$ be the $r$-leaf of $M$ inside $Z$. 
  Recall the Gau\ss{} equation 
  \begin{align*}
    \nabla^{g}_XY=\nabla^{h_r}_X+g(W_r(X),Y)\partial_r, \quad X,Y\in \Gamma(TM)
  \end{align*}
  where $W_r$ denotes the Weingarten map associated to the embedding $\iota_r:M\to M\times \{r\}\subset Z$. 
  For $X\in \Gamma(TM)$ a local coordinate vector field on $M$, we have 
  \begin{align*}
    g(X,\partial_r)=0,\quad [X,\partial_r]=0,\quad  \nabla^{g}_{\partial_r}\partial_r=0,
  \end{align*}
  hence for $m\in M$ are the curves $t\mapsto (t,m)$ geodesics parametrized by the arclength.

  We define the $1$-parameter families of smooth symmetric $2$-tensors $\frac{d}{dr}h_r$, $\frac{d^2}{dr^2}h_r$
  on $M_r$ by
  \begin{align*}
    \frac{d}{dr}h_r(X,Y):=\frac{d}{dr}\big(h_r(X,Y)\big), \quad
    \frac{d^2}{dr^2}h_r(X,Y):=\frac{d^2}{dr^2}\big(h_r(X,Y)\big)
  \end{align*}
  for $X,Y\in T_mM$, which appear in the following set of curvature identities, cf. \cite{BGM}:
  \begin{align*}
   g(W_r(X),Y)=&-\frac{1}{2}\frac{d}{dr}h_r(X,Y),\\
   Ric^{Z}(\partial_r,\partial_r)=&\tr_{h_0}(W_r^2)-\frac{1}{2}\tr_{h_r}(\frac{d^2}{dr^2}h_r),\\
   Ric^{Z}(X,\partial_r)=&X\big(\tr_{h_0}(W_r)\big)-g(\diver^{M}(W_r),X),
  \end{align*}
  for all $X,Y\in T_mM$.

  In fact, our local considerations of a neighborhood of conformal compactification $(X,\bar{g}=dr^2+h_r)$ of  the 
  Poincar\'e-Einstein metric, cf. Appendix B, fit into the concept of the generalized cylinder.
  A tubular neighborhood of the conformal boundary $M$ in $X$ is
  isomorphic to $M\times [0,\varepsilon)$. We can assume that the last coordinate $r\in I$ generates 
  geodesic curve, i.e., the Killing vector frame field $\partial_r$ is parallel with respect to the 
  Levi-Civita connection of $\bar{g}$. 

  For later purposes we collect a few useful facts.
  \begin{lem}\cite[Lemma $6.11.1$, Lemma $6.11.2$]{Juhl}
    
    For $(M,h)$ and the associated Poincar\'e-Einstein metric $g_+$, the two functions $r\mapsto J(dr^2+h_r)$ 
    and $r\mapsto P(dr^2+h_r)$ satisfy
    \begin{align}
      \iota^*\overline{J}=&\iota^* J(dr^2+h_r)=J(h_0)=J,
        \quad \frac{d}{dr}|_{r=0}(\overline{J})=0,\label{eq:RestJ}\\
      \iota^* \overline{P}=&\iota^* P(dr^2+h_r)=P(h_0)=P,\label{eq:RestP}
    \end{align}
    where $\bar{\cdot}$ denotes the evaluation with respect to $\bar{g}=r^2g_+$. 
  \end{lem}

  Let us remind some notation and formulas from Section \ref{Boundary}:  
  \begin{align}
    \widetilde{\slashed D}^{h_r}=&\,
      \partial_r\cdot \sum_{i=1}^n \varepsilon_i s_i\cdot\widetilde{\nabla}^{h_r}_{s_i},\quad
    \overline{\slashed D}=\sum_{i=1}^{n+1} \varepsilon_is_i\cdot\nabla^{\bar{g},S}_{s_i},\notag\\
    \partial_r\cdot \iota_r^*\overline{\slashed D}=&\, 
      \widetilde{\slashed D}^{h_r}\iota_r^*+\frac{n}{2}\iota_r^* H_r-\iota_r^*\nabla^{\bar{g},S}_{\partial_r}
       ,\quad \nabla^{\bar{g},S}_{\partial_r}={\partial_r}+\frac{n}{2}H_r,\label{eq:BicDiracSmallDirac}\\
    H_r=&\frac Jn r-\frac 2n \tr_h(P^2)r^3+\cdots.\notag
  \end{align}
  It is then elementary to verify, using the commutator formula \cite[Proposition $3.1$]{BGM}, 
  the fundamental identity
  \begin{align}
    [\overline{\slashed D},\nabla^{\bar{g},S}_{\partial_r}]=&\frac{1}{2}Ric^X(\partial_r)\cdot
      -\sum_{i=1}^n\varepsilon_i W_r(s_i)\cdot\widetilde{\nabla}^{h_r}_{e_i}-\frac{1}{2}\partial_r\cdot
      \sum_{i=1}^n\varepsilon_i W_r(s_i)\cdot W_r(s_i)\cdot . \label{eq:DiracCov}
  \end{align}
  The formal power series expansions, in the normal coordinate $r$, of the three terms on the 
  right hand side of the last equality are
  \begin{align*}
     \frac{1}{2}Ric^X(\partial_r)\cdot=&\, \frac{1}{2}\sum_{i=1}^n\varepsilon_i Ric^X(\partial_r,s_i)s_i\cdot
       +\frac{1}{2}Ric^X(\partial_r,\partial_r)\partial_r\cdot \\
     =&\, \frac{1}{2}J\partial_r\cdot +O(r^2),\\
     \sum_{i=1}^n\varepsilon_i W_r(s_i)\cdot\widetilde{\nabla}^{h_r}_{s_i}=
       &\, \, r\sum_{i=1}^n\varepsilon_iP(s_i)\cdot\widetilde{\nabla}^{h}_{s_i} +O(r^2), \\
     \frac{1}{2}\partial_r\cdot\sum_{i=1}^n\varepsilon_iW_r(s_i)\cdot W_r(s_i)\cdot =&\, O(r^2),
  \end{align*}
  and collecting all terms together gives
  \begin{align}
     \frac{1}{2}Ric^X(\partial_r)\cdot
       &-\sum_{i=1}^n\varepsilon_iW_r(s_i)\cdot\widetilde{\nabla}^{h_r}_{s_i}
       -\frac{1}{2}\partial_r\cdot\sum_{i=1}^n\varepsilon_iW_r(s_i)\cdot W_r(s_i)\cdot\notag\\
     =&\, \frac{1}{2}J\partial_r\cdot -rh(P,\widetilde{\nabla}^h) + O(r^2).\label{eq:RicciSeries}
  \end{align}

  Based on the review of the fundamental curvature identities for the embedding of 
  the conformal manifold $(M,h)$ into the conformal 
  compactification $(X,\bar{g})$ of the Poincar\'e-Einstein metric $g_+$,
  we are ready to discuss the curved version of factorization identities for 
  residue family operators on spinors of low order. 
  We shall first observe the factorization properties for the first order residue family operator. From the explicit 
  formula in Example \ref{ExamplesCurved} and equation \eqref{eq:BicDiracSmallDirac} we get the following lemma.
  \begin{lem} 
    The residue family operator ${\slashed D}_1^{res}(h; \lambda)$ has the following 
    factorization properties:
    \begin{enumerate}
      \item ${\slashed D}_1^{res}(h; -\frac{n-1}{2})=\widetilde{\slashed D}\iota^*$,
      \item 
        ${\slashed D}_1^{res}(h; -\frac{n}{2})=\partial_r\cdot\iota^*\overline{\slashed D}$.
    \end{enumerate}
    It then follows 
    \begin{align}
      {\slashed D}_1^{res}(h; \lambda)=(2\lambda+n)\widetilde{\slashed D}\iota^*
         -(2\lambda+n-1)\iota^*\partial_r\cdot\overline{\slashed D}.\label{eq:D1}
    \end{align}
  \end{lem} 
  Now we pass to factorization properties for second order residue family operator. 
  \begin{lem} 
    The residue family operator ${\slashed D}_2^{res}(h; \lambda)$ has the following 
    factorization properties:
    \begin{enumerate}
      \item
         ${\slashed D}_2^{res}(h; -\frac{n-3}{2})=\widetilde{\slashed D}{\slashed D}_1^{res}(h; -\frac{n-3}{2})$,
      \item
        ${\slashed D}_2^{res}(h; -\frac{n}{2})=-\partial_r\cdot{\slashed D}_1^{res}(h; -\frac{n+2}{2})
            \overline{\slashed D}$.
    \end{enumerate}
    The residue family operator ${\slashed D}_2^{res}(h; \lambda)$ can be written as 
    \begin{align}
        {\slashed D}_2^{res}(h; \lambda)=(2\lambda+n)\widetilde{\slashed D}^2\iota^* 
             -(2\lambda+n-3)\iota^*\overline{\slashed D}^2
            +2\partial_r\cdot\widetilde{\slashed D}\iota^*\overline{\slashed D}\label{eq:D2}.
    \end{align}
  \end{lem} 
  \begin{bew}
     The first factorization follows from Example \ref{ExamplesCurved}:
     \begin{align*}
        {\slashed D}_2^{res}(h; -\frac{n-3}{2})
            =\widetilde{\slashed{D}}^2\iota^* +2\widetilde{\slashed{D}}\iota^*\partial_r
            =\widetilde{\slashed D}{\slashed D}_1^{res}(h; -\frac{n-3}{2}).
     \end{align*}
     The second one follows from Example \ref{ExamplesCurved} and equation \eqref{eq:D1}:
     \begin{align*}
       {\slashed D}_2^{res}(h; -\frac{n}{2})=&-3\iota^*\nabla^{\bar{g},S}_{\partial_r}
          \nabla^{\bar{g},S}_{\partial_r}
          +\widetilde{\slashed D}^2\iota^*+\frac{3}{2}J\iota^*
          +2\widetilde{\slashed D}\iota^*\nabla^{\bar{g},S}_{\partial_r},\\
       {\slashed D}_1^{res}(h; -\frac{n+2}{2})\overline{\slashed D}
          =&-2\widetilde{\slashed D}\iota^*\overline{\slashed D}
          +3\partial_r\cdot\iota^*\overline{\slashed D}^2\\
       =&2\partial_r\cdot\widetilde{\slashed D}\iota^*\nabla^{\bar{g},S}_{\partial_r}
          -2\partial_r\cdot\widetilde{\slashed D}^2\iota^*+3\partial_r\cdot\widetilde{\slashed D}^2\iota^*\\
       &-3\partial_r\cdot\widetilde{\slashed D}\iota^*\nabla^{\bar{g},S}_{\partial_r}
         -3\iota^*\overline{\slashed D}\nabla^{\bar{g},S}_{\partial_r}
         +3\iota^*[\overline{\slashed D},\nabla^{\bar{g},S}_{\partial_r}].
     \end{align*}
     The combination of equations \eqref{eq:BicDiracSmallDirac}, \eqref{eq:DiracCov} and \eqref{eq:RicciSeries} yields
     \begin{align*}
       {\slashed D}_1^{res}(h; -\frac{n+2}{2})\overline{\slashed D}
         =\partial_r\cdot{\slashed D}_2^{res}(h; -\frac{n}{2}).
     \end{align*}
     As for the explicit formula for ${\slashed D}_2^{res}(h; \lambda)$, it follows from its 
     factorization properties and equation \eqref{eq:D1} that
     \begin{align*}
       {\slashed D}_2^{res}&(h; \lambda)=  \\
       =&\frac{1}{3}(2\lambda+n)\widetilde{\slashed D}{\slashed D}_1^{res}(h; -\frac{n-3}{2})
          +\frac{1}{3}(2\lambda+n-3)\partial_r\cdot{\slashed D}_1^{res}(h; -\frac{n+2}{2})\overline{\slashed D}\\
       =&\frac{1}{3}(2\lambda+n)\widetilde{\slashed D}((-n+3+n)\widetilde{\slashed D}\iota^*-
          (-n+3+n-1)\partial_r\cdot \iota^*\overline{\slashed D})\\
       &+\frac{1}{3}(2\lambda+n-3)\partial_r\cdot((-n-2+n)\widetilde{\slashed D}\iota^*-
          (-n-2+n-1)\partial_r\cdot \iota^*\overline{\slashed D})\overline{\slashed D}\\
       =&(2\lambda+n)\widetilde{\slashed D}^2\iota^*-(2\lambda+n-3)\iota^*\overline{\slashed D}^2
         +2\partial_r\cdot\widetilde{\slashed D}\iota^*\overline{\slashed D}.
     \end{align*}
     The proof is complete. 
  \end{bew}
 
  The factorization properties for the third order residue family operator demonstrate the 
  complexity of the growth of the number of curvature contributions. 
  Before stating the result, we recall the conformal third power of the 
  Dirac operator $\mathcal{D}_3=\slashed{D}^3-2h(P,\nabla^{h,S})-\grad^M(J)\cdot$ on $(M,h)$, 
  cf. references in Appendix A.  
  We denote by $\overline{\mathcal{D}}_3$ the conformal third power of the Dirac operator 
  on $(X,\bar{g})$. Hence, due to the 
  relation of the Clifford multiplication on $S(M,h)$ realized inside $S({X},\bar{g})$, we identify $\mathcal{D}_3$ 
  inside $S(X,\bar{g})$ with $\widetilde{\mathcal{D}}_3$, where all Clifford multiplications are 
  replaced by additional multiplication with 
  $\pm\partial_r\cdot$, cf. Remark \ref{BundleIdentificationCurved}. Then we have   
  \begin{lem}\label{FactThirdOrder}
    The residue family operator ${\slashed D}_3^{res}(h; \lambda)$ on spinors has the following 
    factorization properties:
    \begin{enumerate}
       \item
          ${\slashed D}_3^{res}(h; -\frac{n-3}{2})=\widetilde{{\mathcal D}}_3\iota^*$,
       \item
          ${\slashed D}_3^{res}(h; -\frac{n-5}{2})=\widetilde{\slashed{D}}{\slashed D}_2^{res}(h; -\frac{n-5}{2})$,
       \item
          ${\slashed D}_3^{res}(h; -\frac{n}{2})=\partial_r\cdot{\slashed D}_2^{res}(h; -\frac{n+2}{2})
              \overline{\slashed{D}}$,
       \item
          ${\slashed D}_3^{res}(h; -\frac{n-2}{2})=\partial_r\cdot \iota^*\overline{\mathcal D}_3$.
    \end{enumerate}
  \end{lem} 
  \begin{bew} 
    Let us start with the first identity. From Example \ref{ExamplesCurved} we have
    \begin{align*}
      {\slashed D}_3^{res}(h;-\frac{n-3}{2})=&\widetilde{\slashed D}^3 \iota^*
        -2\partial_r\cdot h(P,\widetilde{\nabla}^h) \iota^*
        +\partial_r\cdot \grad^M(J)\cdot  \iota^*\\
      &-2\partial_r\cdot \grad^M(J)\cdot  \iota^*\\
      =&(\widetilde{\slashed D}^3-2\partial_r\cdot h(P,\widetilde{\nabla}^h)-\partial_r\cdot \grad^M(J)\cdot)\iota^*\\
      =&\widetilde{{\mathcal D}}_3\iota^*.
    \end{align*}
    
    As for the second identity, Example \ref{ExamplesCurved} gives
    \begin{align*}
      \widetilde{\slashed{D}}{\slashed D}_2^{res}(h; -\frac{n-5}{2})=&
         \widetilde{\slashed D}(2\iota^*(\nabla^{\bar{g},S}_{\partial_r})^2-\widetilde{\slashed D}^2\iota^*-J\iota^*
               +2\widetilde{\slashed D}\iota^*\nabla^{\bar{g},S}_{\partial_r})\\
      =& 2\widetilde{\slashed D}\iota^*(\nabla^{\bar{g},S}_{\partial_r})^2
               +\widetilde{\slashed D}^3\iota^*-\partial_r\cdot \grad^M(J)\cdot \iota^*\\
      &-J\widetilde{\slashed D}\iota^*+2\widetilde{\slashed D}^2\iota^*\nabla^{\bar{g},S}_{\partial_r}\\
       =&{\slashed D}_3^{res}(h; -\frac{n-5}{2}).
    \end{align*}
    
    To prove the third identity, we need Example \ref{ExamplesCurved} again and equations \eqref{eq:BicDiracSmallDirac}, 
    \eqref{eq:DiracCov}, \eqref{eq:RicciSeries}. On the one hand,
    \begin{align*}
      {\slashed D}_3^{res}(h; -\frac{n}{2})=&5\iota^*(\nabla^{\bar{g},S}_{\partial_r})^3
         -3\widetilde{\slashed D}^2\iota^*\nabla^{\bar{g},S}_{\partial_r}
         -\frac{15}{2}J\iota^*\nabla^{\bar{g},S}_{\partial_r}
         -3\widetilde{\slashed D}\iota^*(\nabla^{\bar{g},S}_{\partial_r})^2\\
       &+\widetilde{\slashed D}^3\iota^*-5\partial_r\cdot h(P,\widetilde{\nabla}^h)\iota^*
         +\frac{3}{2}J\widetilde{\slashed D}\iota^*-\partial_r\cdot \grad^M(J)\cdot \iota^*,
     \end{align*}
     while
     \begin{align*} 
       {\slashed D}_2^{res}&(h; -\frac{n+2}{2})\overline{\slashed D}=
          (-5\iota^*(\nabla^{\bar{g},S}_{\partial_r})^2 + \widetilde{\slashed D}^2\iota^* +\frac{5}{2}J\iota^*
          +2\widetilde{\slashed D}\iota^*\nabla^{\bar{g},S}_{\partial_r})\overline{\slashed D}\\
       =&-5\iota^*(\nabla^{\bar{g},S}_{\partial_r})^2\overline{\slashed D}+
          \widetilde{\slashed D}^2(-\partial_r\cdot\widetilde{\slashed D}\iota^*
          +\partial_r\cdot \iota^*\nabla^{\bar{g},S}_{\partial_r})
          +2\widetilde{\slashed D}\iota^*\nabla^{\bar{g},S}_{\partial_r}\overline{\slashed D}\\
        &+\frac{5}{2}J(-\partial_r\cdot\widetilde{\slashed D}\iota^*
          +\partial_r\cdot \iota^*\nabla^{\bar{g},S}_{\partial_r})\\          
        =&5\partial_r\cdot\widetilde{\slashed D}\iota^*(\nabla^{\bar{g},S}_{\partial_r})^2
               -5\partial_r\cdot \iota^*(\nabla^{\bar{g},S}_{\partial_r})^3 
               +5\iota^* [\overline{\slashed D},\nabla^{\bar{g},S}_{\partial_r}]\nabla^{\bar{g},S}_{\partial_r}\\
         &+5\iota^* \nabla^{\bar{g},S}_{\partial_r}[\overline{\slashed D},\nabla^{\bar{g},S}_{\partial_r}]
                -\partial_r\cdot\widetilde{\slashed D}^3\iota^*
                +\partial_r\cdot\widetilde{\slashed D}^2\iota^*\nabla^{\bar{g},S}_{\partial_r}
                -\frac{5}{2}J\partial_r\cdot\widetilde{\slashed D}\iota^* \\
          &+\frac{5}{2}J\partial_r\cdot \iota^*\nabla^{\bar{g},S}_{\partial_r}
                +2\partial_r\cdot\widetilde{\slashed D}^2\iota^*\nabla^{\bar{g},S}_{\partial_r}
                -2\partial_r\cdot\widetilde{\slashed D}\iota^*(\nabla^{\bar{g},S}_{\partial_r})^2
                -2\widetilde{\slashed D}\iota^* [\overline{\slashed D},\nabla^{\bar{g},S}_{\partial_r}]\\
           =&3\partial_r\cdot\widetilde{\slashed D}\iota^*(\nabla^{\bar{g},S}_{\partial_r})^2
                -5\partial_r\cdot \iota^*(\nabla^{\bar{g},S}_{\partial_r})^3
                +3\partial_r\cdot\widetilde{\slashed D}^2\iota^*\nabla^{\bar{g},S}_{\partial_r}
                 -\partial_r\cdot\widetilde{\slashed D}^3\iota^*\\
            &-\frac{3}{2}J\partial_r\cdot\widetilde{\slashed D}\iota^*-\grad^M(J)\cdot \iota^*
                +5J\partial_r\cdot \iota^*\nabla^{\bar{g},S}_{\partial_r}
                +5\iota^*\nabla^{\bar{g},S}_{\partial_r}[\overline{\slashed D},\nabla^{\bar{g},S}_{\partial_r}].
       \end{align*}
       After multiplication by $\partial_r\cdot$, we arrive at
       \begin{align*}
          {\slashed D}_3^{res}(h; -\frac{n}{2})=\partial_r\cdot
              {\slashed D}_2^{res}(h; -\frac{n+2}{2})\overline{\slashed D}.
       \end{align*}

       Finally we come to the last identity. Again, Example \ref{ExamplesCurved} implies
       \begin{align*}
         {\slashed D}_3^{res}(h; -\frac{n-2}{2})=& \iota^*(\nabla^{\bar{g},S}_{\partial_r})^3
           -\widetilde{\slashed D}^2\iota^*\nabla^{\bar{g},S}_{\partial_r}
           -\frac{3}{2}J\iota^*\nabla^{\bar{g},S}_{\partial_r}
                -\widetilde{\slashed D}\iota^*(\nabla^{\bar{g},S}_{\partial_r})^2+\widetilde{\slashed D}^3\iota^*\\
         & -3\partial_r\cdot h(P,\widetilde{\nabla}^h)\iota^* 
            +\frac{1}{2}J\widetilde{\slashed D}\iota^*
                -\partial_r\cdot \grad^M(J)\cdot \iota^*.
       \end{align*}
       Because
       \begin{align*}
             \overline{{\mathcal D}}_3=\overline{\slashed D}^3-2\bar{g}(\overline{P},\nabla^{\bar{g},S})
                -\grad^{X}(\overline{J})\cdot,
       \end{align*}
       we have
       \begin{align*}
           \iota^*\overline{{\mathcal D}}_3=\iota^*\overline{\slashed D}^3
            -2\iota^*\bar{g}(\overline{P},\nabla^{\bar{g},S})-\iota^* \grad^{X}(\overline{J}).
       \end{align*}
       From equations \eqref{eq:RestJ} and \eqref{eq:RestP} we conclude
       \begin{align*}
         -2\iota^*\bar{g}(\overline{P},\nabla^{\bar{g},S})
                  =&-2\iota^*\sum_{i=1}^{n+1}\overline{P}(s_i)\cdot \nabla^{\bar{g},S}_{s_i}
                  =-2{h}({P},\widetilde{\nabla}^h)\iota^*,\\
          -\iota^* \grad^{X}(\overline{J})\cdot
                  =&-\iota^* (\sum_{i=1}^{n}s_i(\overline{J})s_i\cdot +\partial_r(\overline{J})\partial_r\cdot)
                  =-\grad^{M}({J})\iota^*\cdot.
       \end{align*}
       Furthermore, using equations \eqref{eq:BicDiracSmallDirac}, \eqref{eq:DiracCov} and 
       \eqref{eq:RicciSeries} we compute
       \begin{align*}
          \iota^*\overline{{\slashed D}}_3=&-\partial_r\cdot\widetilde{\slashed D}
                (-\partial_r\cdot\widetilde{\slashed D}\iota^*
                +\partial_r\cdot \iota^*\nabla^{\bar{g},S}_{\partial_r})\overline{\slashed D}
                +\partial_r\cdot \iota^*\nabla^{\bar{g},S}_{\partial_r}\overline{\slashed D}^2\\
           =&-\partial_r\cdot\widetilde{\slashed D}^3\iota^* 
                + \partial_r\cdot\widetilde{\slashed D}^2\iota^*\nabla^{\bar{g},S}_{\partial_r}
                 -\partial_r\cdot\widetilde{\slashed D}^2\iota^*\nabla^{\bar{g},S}_{\partial_r}
                +\partial_r\cdot\widetilde{\slashed D}\iota^*(\nabla^{\bar{g},S}_{\partial_r})^2\\
            &+\widetilde{\slashed D}\iota^* [\overline{\slashed D},\nabla^{\bar{g},S}_{\partial_r}]
                +\partial_r\cdot (-\partial_r\cdot\widetilde{\slashed D}\iota^* 
                +\partial_r\cdot \iota^*\nabla^{\bar{g},S}_{\partial_r})
                \overline{\slashed D}\nabla^{\bar{g},S}_{\partial_r}\\
             &-\partial_r\cdot(-\partial_r\cdot\widetilde{\slashed D}\iota^* 
                +\partial_r\cdot \iota^*\nabla^{\bar{g},S}_{\partial_r})
                [\overline{\slashed D},\nabla^{\bar{g},S}_{\partial_r}]
                -\partial_r\cdot \iota^* [\overline{\slashed D},\nabla^{\bar{g},S}_{\partial_r}]\overline{\slashed D}\\
             =& -\partial_r\cdot\widetilde{\slashed D}^3\iota^*
                + \partial_r\cdot\widetilde{\slashed D}^2\iota^*\nabla^{\bar{g},S}_{\partial_r}
                +\partial_r\cdot\widetilde{\slashed D}\iota^*(\nabla^{\bar{g},S}_{\partial_r})^2
                 -\partial_r\cdot \iota^*(\nabla^{\bar{g},S}_{\partial_r})^3\\
             &+J\partial_r\cdot \iota^*\nabla^{\bar{g},S}_{\partial_r}
                -\frac{1}{2}J\partial_r\cdot\widetilde{\slashed D}\iota^* 
                +\frac{1}{2}J\partial_r\cdot \iota^* \nabla^{\bar{g},S}_{\partial_r}
                -h(P,\widetilde{\nabla}^h)\iota^*.
       \end{align*}
       Taking into account all contributions we get
       \begin{align*}
          \iota^*\overline{{\mathcal D}}_3=&-\partial_r\cdot\widetilde{\slashed D}^3\iota^*
                +\partial_r\cdot\widetilde{\slashed D}^2\iota^*\nabla^{\bar{g},S}_{\partial_r}
                +\partial_r\cdot\widetilde{\slashed D}\iota^*(\nabla^{\bar{g},S}_{\partial_r})^2
                -\partial_r\cdot \iota^*(\nabla^{\bar{g},S}_{\partial_r})^3\\
          &+\frac{3}{2}J\partial_r\cdot \iota^*\nabla^{\bar{g},S}_{\partial_r} 
                -\frac{1}{2}J\partial_r\cdot \widetilde{\slashed D} \iota^*-3h({P},\widetilde{\nabla}^h)\iota^* 
                -\grad^{M}({J})\cdot \iota^* 
       \end{align*}
       and the proof is complete.
  \end{bew}
  
  From the Lemma above and Theorem \ref{ConformalTrafoLawResidueFamily} we conclude:
  \begin{kor}
    Let $(M,h)$ be a semi-Riemannian $Spin$-manifold. 
    Then the conformal first and third power of the Dirac operator are given by 
    \begin{align*}
       \mathcal{D}_1=&\slashed{D},\\
       \mathcal{D}_3=&\slashed{D}^3-2h(P,\nabla^{h,S})-\grad^M(J)\cdot.
    \end{align*}
  \end{kor}
  \begin{bew}
    The results follow from the factorizations of $\slashed{D}^{res}_1(h;-\frac{n-1}{2})$ and 
    $\slashed{D}^{res}_3(h;-\frac{n-3}{2})$ by post- and pre-composition with the identification map $\Xi$ 
    and its inverse, cf. Remark \ref{BundleIdentificationCurved}. 
    Note that one could alternatively use $\Xi:=\Xi^-$ in the case of odd $n$, with the effect of an additional 
    sign.  
  \end{bew}

  \begin{bem}
    The third order residue family operator is a polynomial in $\lambda$ of degree two, cf.  
    Example \ref{ExamplesCurved}. Lemma \ref{FactThirdOrder} 
    shows that it satisfies four factorization identities. Thus making 
    an ansatz for $\slashed{D}^{res}_3(h;\lambda)$ as a polynomial of degree three in $\lambda$ 
    shows that due to four factorizations the operator valued coefficient by the third power of 
		$\lambda$ is trivial. However, this operator 
    is a multiple of 
    \begin{align*}
      M_3\iota^*-\iota^*\overline{M}_3,
    \end{align*}
    where $M_3:=-2h(P,\nabla^{h,S})-\grad^M(J)\cdot$ and $\overline{M}_3$ is given by analogous 
    formula evaluated with respect to $\bar{g}$. 
    Notice that this operator is one of the three operators $\{M_1,M_3,M_5\}$ 
    described in \cite[Chapter $6$]{Fischmann}, which determine the conformal powers of the 
    Dirac operator, up to order five, as a non-commutative free algebra.  
  \end{bem}

  Now we discuss the factorization identities for $D_1(X,M;g,\lambda)$, see \eqref{eq:FirstOrderFamilyGeneralHyper} for 
  its definition. Based on equation \eqref{eq:BicDiracSmallDirac} and direct computations we obtain 
  \begin{lem}
    The family of first order operators $D_1(X,M;g,\lambda)$ has the following factorization properties:
    \begin{enumerate}      
        \item $D_1(X,M;g,-\frac{n-1}{2})=N(g)\cdot\iota^*\slashed{D}^g$, 
        \item $D_1(X,M;g,-\frac{n}{2})=\widetilde{\slashed{D}}^h\iota^*$. 
    \end{enumerate}
    The family $D_1(X,M;g,\lambda)$ can be written as
    \begin{align*}
      D_1(X,M;g,\lambda)=&-(2\lambda+n-1)N(g)\cdot\iota^*\slashed{D}^g
        +(2\lambda+n)\widetilde{\slashed{D}}^h\iota^*\\
      &+2(\lambda+\frac n2)(\lambda+\frac{n-1}{2})\iota^* H(g).
    \end{align*}
  \end{lem}


\section{Representation theory and residue family operators}\label{SingularVectors}

  In the present section we review the interpretation of our results, concerning factorization identities, 
  in the framework of the classification 
  of homomorphisms of generalized Verma modules for couples of Lie algebras 
  ${\mathfrak g}$, ${\mathfrak g}'$ and their conformal parabolic subalgebras ${\mathfrak p}$,
  ${\mathfrak p}'$. As for the technique and general background we refer to \cite{KOSS}. 
  The singular vectors describing such homomorphisms
  correspond to conformally covariant differential operators in the pairing of generalized Verma 
  modules with induced representations.

  Let ${\mathfrak g}$ resp. ${\mathfrak g}'$ be the conformal Lie algebras in dimension $n$
  resp. $n-1$.
  Let $\lambda \in {\mathbb{C}}$, $N \in {\mathbb{N}}_0$, and let $s_\lambda\in{\mathbb S}_\lambda$ 
  for the spinor representation ${\mathbb S}_\lambda=\Delta_{n-1}\otimes{\mathbb C}_\lambda$ of 
  $\R^{n-1}$ tensored by $1$-dimensional representation 
  of the Levi factor of ${\mathfrak p}$ on $\C_\lambda$, $a\mapsto a^{\lambda}, \, a\in{\mathbb R}$. 
  It is a result of \cite{KOSS} that ${\mathfrak g}'$-singular vectors of odd and even homogeneity are given 
  by two Gegenbauer polynomials, cf. Appendix C:
  \begin{align}
    \tilde{P}_N(t):=&(-t)^NC^{-\lambda-\frac n2}_{2N}\left(\frac{i}{\sqrt{t}}\right)\notag\\
    =&\frac{1}{N!}\left(-(\lambda+\frac{1}{2})-\frac{n-1}{2}\right)_N
    \sum_{j=0}^N a_j^{(N)}(\lambda+\frac{1}{2})t^j,\label{eqn:5.7}\\
    \tilde{Q}_{N}(t):=&i(-t)^N\frac{1}{\sqrt{t}}C^{-\lambda-\frac n2}_{2N+1}\left(\frac{i}{\sqrt{t}}\right)\notag\\
    =&\frac{2}{N!}\left(-(\lambda+\frac{1}{2})-\frac{n-1}{2}\right)_{N+1}
      \sum_{j=0}^{N}b_j^{(N)}(\lambda+\frac{1}{2})t^j.\label{eqn:5.8}
  \end{align} 
  In particular, the polynomial in $t$ (and some other variables $\eta,{\underline\xi},{\underline \eta}$) 
  valued in the algebra of endomorphisms of ${\mathbb S}_\lambda$,
  \begin{align*}
     \tilde{F}_{2N}\cdot s_\lambda  = (\eta^{2N}\tilde{P}_N(t)+\tilde{Q}_{N-1}(t)
       {\underline\xi}\cdot{\underline \eta})\cdot s_\lambda ,
  \end{align*}
  determines a ${\mathfrak g}'$-singular vector of homogeneity $2N$, 
  $N\in{\mathbb N}_0$, in the Fourier dual of the 
  generalized Verma ${\mathfrak g}$-module induced from the twisted spinor representation of the conformal parabolic 
  subalgebra ${\mathfrak p}$. Note that $\tilde{Q}_{-1}(t):=0$ by convention. 
  Furthermore, it follows from the system of differential equations satisfied by 
  $\tilde{P}_N(t)$ and $\tilde{Q}_{N-1}(t)$, cf. \cite{KOSS}, that 
  the coefficients of $\tilde{P}_N(t),\tilde{Q}_{N-1}(t)$ are related by
  \begin{align}
      N(2j-2N+1)b_j^{(N-1)}(\lambda+\frac{1}{2})-(j+1)a_{j+1}^{(N)}(\lambda+\frac{1}{2})=0\label{eq:EvenNorm1}
  \end{align}
  for all $j=0,\dots ,N-1$, and
  \begin{align}
      N(2\lambda+n-4N+2j+2)b_j^{(N-1)}(\lambda+\frac{1}{2})
       +(j-N)a_{j}^{(N)}(\lambda+\frac{1}{2})=0\label{eq:EvenNorm2}
  \end{align}
  for all $j=0,\dots ,N$.
 
  If we did not fix the normalizations of Gegenbauer polynomials, cf. Subsection \ref{ResFamFlatCase}, 
  the relation between the two sets of coefficients would be $a^{(N)}_N(\lambda)=-b_{N-1}^{(N-1)}(\lambda)$.   

  In the case of odd homogeneity $2N+1$, $N\in{\mathbb N}_0$,
    \begin{align}
     \tilde{F}_{2N+1}\cdot s_\lambda  = \eta^{2N}(\tilde{P}_N(t){\underline \xi}+\tilde{Q}_N(t)
       {\underline \eta} )\cdot s_\lambda
  \end{align}
  determine ${\mathfrak g}'$-singular vectors provided the coefficients of $\tilde{P}_N(t), \tilde{Q}_N(t)$ 
  satisfy
    \begin{align}
      -(2\lambda+n-2N)(2N+1-2j)&b_j^{(N)}(\lambda+\frac{1}{2})\nonumber \\
        &+(2\lambda +n-4N+2j)a_{j}^{(N)}(\lambda+\frac{1}{2})=0\label{eq:OddNorm1}
    \end{align}
    for all $j=0,\dots ,N$, and 
    \begin{align}
        (N-j)a_j^{(N)}(\lambda+\frac{1}{2})
         +(j+1)(2\lambda+n-2N)b_{j+1}^{(N)}(\lambda+\frac{1}{2})=0\label{eq:OddNorm2}
    \end{align}
    for all $j=0,\dots ,N-1$.

    If we did not fix the normalizations of Gegenbauer polynomials, cf. Subsection \ref{ResFamFlatCase}, 
    the relation between the two sets of coefficients would be $a^{(N)}_N(\lambda)=b_{N}^{(N)}(\lambda)$.   

  \begin{bem}
    The relations among the coefficients of Gegenbauer polynomials 
    can be used to prove the factorization identities of residue family operators stated in 
    Theorem \ref{FactorziationBicDirac}. For example, the left hand side of both 
    equations \eqref{eq:check2} and \eqref{eq:check6} satisfies a recurrence relation, cf. equation \eqref{eq:a-recurrence} 
    and \eqref{eq:b-recurrence}. 
    The right hand side of both equations \eqref{eq:check2} and \eqref{eq:check6} 
    satisfies the same recurrence relation due to equations 
    \eqref{eq:EvenNorm1}, \eqref{eq:EvenNorm2}, \eqref{eq:OddNorm1} and \eqref{eq:OddNorm2}. 
  \end{bem}


\section{Poisson transformation and residue family operators}\label{Poisson}

In the present section we indicate the origin of the family of distributions 
\begin{align*}
  \delta_N(h;\lambda):\mathcal{C}_c^\infty(\R^n_{\geq 0})\to \mathcal{C}^\infty(\R^{n-1})
\end{align*}
on the real hyperbolic space (upper half-space), whose residues produce the 
residue families on functions, cf. \cite{Juhl}. We also comment on the relation between $\delta_N(h;\lambda)$ 
and Poisson transformation acting on induced representations on the boundary of the rank one 
symmetric space 
\begin{align*}
  \Ham^n:=\big(\R^n_{x_n>0},g_{hyp}:=x_n^{-2}(dx_n^2+dx_1^2+\cdots+dx_{n-1}^2)\big).
\end{align*}
All representation spaces are considered in the non-compact picture. In the second part of the 
section we pass to the spinor valued case and suggest analogous distribution valued in the algebra
of endomorphisms of spinor representation.

We shall start with the scalar case, so we denote by $x$ a point in the vector space $\R^n$
and write it as $x=(x^\prime,x_n)$ with 
respect to the splitting $\R^n=\R^{n-1}\times\R$. 
The Poisson kernel on $\Ham^n$ is given by locally integrable function
\begin{align*}
  P(y,x^\prime):=\frac{y_n}{\abs{(y^\prime,y_n)-(x^\prime,0)}^2},\quad y\in\R^n, x^\prime\in\R^{n-1}, 
\end{align*}
where the absolute value is considered with respect to the euclidean scalar product. The 
Poisson transform $\mathcal{P}_\mu:\mathcal{C}_c^\infty(\R^{n-1})\to{\mathcal C}^\infty(\R^n)$
with Poisson kernel $P(y,x^\prime)$, given by
\begin{align*}
  \mathcal{P}_\mu(f)(y):=\int_{\R^{n-1}} P(y,x^\prime)^\mu f(x^\prime)dx^\prime, \quad\mu\in{\mathbb C},
\end{align*}
is an eigenfunction of the hyperbolic Laplace operator, see \cite[Theorem $1.7$]{Helgason},
\begin{align*}
-\Delta_{g_{hyp}}u=\mu(n-1-\mu)u,
\end{align*}
for $u\in\mathcal{C}^\infty(\R^n)$. 

For $f\in\mathcal{C}_c^\infty(\R^{n-1})$, we denote by $u:=\mathcal{P}_\mu(f)$ its Poisson transform.
Let us consider the family of locally integrable functions $x_n^{\lambda-n}u$ for $Re(\lambda)\gg 0$,
and define the family of distributions $M_u(\lambda;x_n)$ supported on $\R^n_{\geq 0}$ by 
\begin{align*}
  \big(M_u(\lambda;y_n)\big)(\varphi)
   =\int_{\R^n_{\geq 0}}y_n^{\lambda-n}u(y^\prime,y_n)\varphi(y^\prime,y_n)dy^\prime dy_n
\end{align*}
for compactly supported $\varphi\in\mathcal{C}^\infty_c(\R^n_{\geq 0})$, $Re(\lambda)\gg 0$. 
Because $u$ is the Poisson transform of 
$f\in\mathcal{C}_c^\infty(\R^{n-1})$, we get
\begin{align}
  \big(M_u(\lambda;y_n)\big)(\varphi)
    =\int_{\R^{n-1}}\int_{\R^n_{\geq 0}}\varphi(y^\prime,y_n)
    \frac{y_n^{\lambda-n+\mu}}{(\abs{y^\prime-x^\prime}^2+y_n^2)^\mu}f(x^\prime)
    dy^\prime dy_n dx^\prime.\label{eq:PoisonDiscription}
\end{align}
As noticed in \cite{kob}, for $\lambda,\mu\in \C$ such that $Re(\lambda-\mu)>0$ and 
$Re(\lambda+\mu)>n-1$, the locally integrable function 
\begin{align}
  K_{\lambda,\mu}(x^\prime,x_n)
      :=\frac{\abs{x_n}^{\lambda+\mu-n}}{(\abs{x^\prime}^2+x_n^2)^\mu}\label{scdistrfam}
\end{align}
allows to introduce an integral operator 
\begin{align*}
  \mathcal{K}_{\lambda,\mu}:{\mathcal C}_c^\infty(\R^n)&\to{\mathcal C}^\infty(\R^{n-1})\\
  f&\mapsto (\mathcal{K}_{\lambda,\mu}f)(x^\prime)
    := \int_{\R^n} f(y^\prime,y_n) K_{\lambda,\mu}(x^\prime-y^\prime,-y_n) dy^\prime dy_n 
\end{align*}
intertwining the action of the conformal group associated to the boundary $\R^{n-1}$. 
Regarding the image of that integral operator as a distribution supported on $\R^{n-1}$, we can 
evaluate it on $g\in {\mathcal C}^\infty_c(\R^{n-1})$,
\begin{align}
   (\mathcal{K}_{\lambda,\mu}f)(g)
     =&\int_{\R^{n-1}}(\mathcal{K}_{\lambda,\mu}f)(x^\prime) g(x^\prime) dx^\prime\notag\\
   =&\int_{\R^{n-1}}\int_{\R^n}f(y^\prime,y_n)K_{\lambda,\mu}(x^\prime-y^\prime,-y_n) g(x^\prime) 
     dy^\prime dy_n dx^\prime.\label{eq:KernelDiscription}
\end{align} 
A direct comparison of equations \eqref{eq:PoisonDiscription} and \eqref{eq:KernelDiscription} gives 
\begin{sat}\label{ResidueKernel}
  Let $f\in{\mathcal C}_c^\infty(\R^{n-1})$, $\varphi\in\mathcal{C}^\infty_c(\R^n_{\geq 0})$ be
  compactly supported functions, $u=\mathcal{P}_\mu(f)$. Then
  \begin{align}\label{distgivenbyfamilykernel}
    \big(M_u(\lambda,x_n)\big)(\varphi)=(\mathcal{K}_{\lambda,\mu}\varphi)(f),
  \end{align}
  and the equality extends continuously to the space of smooth functions.
\end{sat}
The meromorphic continuation of the distribution $\mathcal{K}_{\lambda,\mu}$ to $\lambda,\mu\in{\mathbb C}$,
\cite{kob}, yields the residue families on densities, cf. \cite{Juhl}. 

Now we turn our attention to an analogous question on Poisson transformation for spinor valued functions, i.e., 
we suggest a family of distributions whose meromorphic continuation presumably yields the residue family operators 
on spinors discussed in the present article. 

The Poisson kernel for spinors on the real hyperbolic space ${\mathbb H}^n$ is a locally integrable function
valued in the Clifford algebra $Cl({\mathbb R}^n)$,
\begin{eqnarray}\label{pks}
  P^S_\mu(y,x^\prime):=\frac{y_n^\mu}{\abs{(y^\prime,y_n)-(x^\prime,0)}^{{n}}}
  (\sum_{i=1}^{n-1}(y_i-x_i)e_i+y_ne_n)e_n\cdot ,\, \mu\in{\mathbb C},
\end{eqnarray}
cf. \cite{cam}.
Here $y\in\R^n, x^\prime\in\R^{n-1}$, where the absolute value is given by the 
euclidean scalar product. 
\begin{sat}
  The Poisson transformation 
  \begin{align*}
   \mathcal{P}^S_\mu:\Gamma_c\big(S({\mathbb R}^{n},h)|_{\R^{n-1}}\big)\to\Gamma\big(S(\R^{n},\bar{g})\big)
  \end{align*}
  with Poisson kernel $P^S_\mu(y,x^\prime)$, given by
  \begin{align*}
    \mathcal{P}^S_\mu(s)(y):=\int_{\R^{n-1}} P_\mu^S(y,x^\prime) s(x^\prime)dx^\prime,
  \end{align*}
  is an eigenvector of the hyperbolic Dirac operator,
  \begin{align*}
    {\slashed D}^{g_{hyp}}\mathcal{P}^S_\mu(s)=(\lambda -\frac{n-1}{2})\mathcal{P}^S_\mu(s), \,\, 
    s\in \Gamma_c\big(S(\R^{n},h)|_{\R^{n-1}}\big),
  \end{align*}
  with the eigenvalue $\lambda -\frac{n-1}{2}$.
\end{sat}
\begin{bew}
  The proof is a direct consequence of the explicit formula
  \begin{align*}
     \slashed{D}^{g_{hyp}}=x_ne_n\cdot\partial_n-\frac 12(n-1)e_n\cdot+x_n\sum_{i=1}^{n-1}e_i\cdot\partial_i,
  \end{align*}
  the hyperbolic Dirac operator, and the definition of $\mathcal{P}^S_\mu(s)$.
\end{bew}

Let us denote by $D\big({S}({\mathbb H}^n)\big)$ the algebra of invariant differential operators
acting on spinors for rank one symmetric space ${\mathbb H}^n$. It is well-known, cf. \cite{cam}, that for $n$ even, 
$D\big({S}^\pm({\mathbb H}^n)\big)\simeq\C[({{\slashed D}^{g_{hyp}}})^2]|_{{S}^\pm({\mathbb H}^n)}$,
and for $n$ odd, ${D}\big({S}({\mathbb H}^n)\big)\simeq\C[{\slashed D}^{g_{hyp}}]$. 
Consequently, the Poisson 
transform maps smooth spinors on the boundary of ${\mathbb H}^n$ into eigenspaces of the Dirac 
operator on the hyperbolic space ${\mathbb H}^n$.

The distribution on spinors \eqref{pks} allows a generalization along the lines of the scalar 
valued distribution, cf. equation \eqref{scdistrfam}. 
We define, for $\lambda,\mu\in \C$ such that $Re(\lambda)\gg 0$ and 
$Re(\mu)\gg 0$, a locally integrable $Cl({\mathbb R}^n)$-valued function on ${\R}^n\times\R^{n-1}$
\begin{align}
  K^S_{\lambda,\mu}(x^\prime,x_n):=\frac{\abs{x_n}^{\lambda}}{(\abs{x^\prime}^2+x_n^2)^{\frac{n}{2}}}
  \big((\sum_{i=1}^{n-1}(y_i-x_i)e_i+y_ne_n)e_n\cdot\big)^{\mu},\label{spinorcdistrfam}
\end{align}
which allows to introduce an integral operator on spinors 
\begin{align*}
  \mathcal{K}^S_{\lambda,\mu}:\Gamma_c\big( S(\R^n,\bar{g}) \big)&\to\Gamma\big(S(\R^{n},h)|_{\R^{n-1}}\big),\\
  s & \mapsto (\mathcal{K}^S_{\lambda,\mu}\cdot s)(x^\prime)
    := \int_{\R^n} K_{\lambda,\mu}(x^\prime-y^\prime,-y_n)\cdot s(y^\prime,y_n)dy^\prime dy_n. 
\end{align*}
Regarding the image of the previous integral operator as a distribution supported on $\R^{n-1}$, 
the integration of a smooth compactly supported spinor valued function 
$\varphi\in \Gamma_c\big(S(\R^{n},h)|_{\R^{n-1}}\big)$ yields
\begin{align}
   (\mathcal{K}^S_{\lambda,\mu}s)(\varphi)
     =&\int_{\R^{n-1}}<\mathcal{K}^S_{\lambda,\mu}\cdot s(x^\prime), \varphi(x^\prime)> dx^\prime\notag\\
   =&\int_{\R^{n-1}}\int_{\R^n}<K^S_{\lambda,\mu}(x^\prime-y^\prime,-y_n)\cdot s(y^\prime,y_n),
      \varphi(x^\prime)>dy^\prime dy_n dx^\prime.\label{eq:KernelDiscriptionspinor}
\end{align} 
It remains to explain the definition of the family of distributions 
\begin{align*}
  \big((\sum_{i=1}^{n-1}(y_i-x_i)e_i+y_ne_n)e_n\cdot\big)^{\mu},\quad \mu\in{\mathbb C}.
\end{align*}
Let us define two $Cl({\mathbb R}^n)$-valued functions 
$\Pi_\pm:=\frac{1}{2}(Id\pm i\omega)$,
$\omega\in S^{n-1}\subset{\mathbb R}^n\subset Cl({\mathbb R}^n)$, fulfilling $-i\omega\Pi_\pm=\pm\Pi_\pm$, or 
\begin{align*}
  {(\sum_{i=1}^{n-1}(y_i-x_i)e_i+y_ne_n)e_n}=r\frac{{(\sum_{i=1}^{n-1}(y_i-x_i)e_i+y_ne_n)e_n}}
     {\abs{(\sum_{i=1}^{n-1}(y_i-x_i)e_i+y_ne_n)e_n}}
\end{align*}
for $\omega:=\frac{{(\sum_{i=1}^{n-1}(y_i-x_i)e_i+y_ne_n)e_n}}{\abs{(\sum_{i=1}^{n-1}(y_i-x_i)e_i+y_ne_n)e_n}}\in S^{n-1}$
in the spherical coordinates on $\R^n$ with the radial coordinate $r$. 
Then 
\begin{align*}
  (-i{(\sum_{i=1}^{n-1}(y_i-x_i)e_i+y_ne_n)e_n})^\mu\Pi_\pm
     =(-i\omega)^\mu r^\mu\Pi_\pm=(\pm 1)^\mu r^\mu\Pi_\pm,
\end{align*}
and because $\Pi_++\Pi_-=Id$, we get 
\begin{align*}
  (-i(\sum_{i=1}^{n-1}&(y_i-x_i)e_i+y_ne_n)e_n)^\mu=1^\mu r^\mu\Pi_++(-1)^\mu r^\mu\Pi_- \\
  =&\frac{1}{2}(1^\mu +(-1)^\mu)r^\mu +\frac{1}{2}(1^\mu -(-1)^\mu)(-i{(\sum_{i=1}^{n-1}(y_i-x_i)e_i+y_ne_n)e_n})r^{\mu-1} .
\end{align*}
As for the distribution $r^\mu$, we have the standard definition as in the scalar valued case,
\begin{align*}
  (r_+^\mu,\varphi)=\int_{0}^\infty r_+^\mu\varphi(r)dr,\quad \varphi\in {\mathcal C}^\infty_c({\mathbb R}_+),
\end{align*}
locally integrable for $Re(\mu)>-1$, and 
\begin{align*}
   (1)^\mu=e^{2\pi il\mu},\, (-1)^\mu=e^{(2l+1)\pi i\mu},\quad l\in{\mathbb Z}.
\end{align*}
In this way, \eqref{spinorcdistrfam} reduces to a linear combination of two scalar valued distributions
with coefficients in $Cl({\mathbb R}^n)$.

We expect that this distribution will play a decisive role in the spinor version of the scalar distribution kernel, 
cf. Proposition \ref{ResidueKernel}. 

\section{Discussion and outlook}\label{Outlook}

In this last short section we comment on several unresolved questions related to our work. 

For a manifold with conformal structure $(M,[h])$, there is an invariant theory of tractors based
on the existence of conformally invariant connection and compatible invariant 
metric on indecomposable bundles called tractor bundles. This structure allows, for a hypersurface $M$ inside $(X,g)$, 
to produce natural
intrinsically conformally covariant operators for $N\in{\mathbb N}_0$ and $\lambda\in{\mathbb C}$,
\begin{align*}
  D_{N}^{T}(X,M;g,\lambda):\mathcal{C}^\infty(X)\to\mathcal{C}^\infty(M)
\end{align*}
(the superscript $T$ stands for the tractor construction), cf. \cite[Definitions $6.21.1$ and $6.21.2$]{Juhl}. 
This family specializes, in the case of conformal compactification $(X,\bar{g})$ of the Poincar\'e-Einstein metric $g_+$ 
associated to conformally flat manifold $(M,h)$, to residue families:
\begin{align*}
  D_{2N}^T(X,M;\bar{g},\lambda)=D^{res}_{2N}(h; \lambda),
\end{align*}
cf. \cite[Theorem $6.21.2$]{Juhl}. We conjecture that an analogous statement holds for conformally covariant 
operators acting on sections of a spinor bundle.  

Another interesting point is the question on the origin of the spinor residue family operators.
In \cite{kob}, the Riesz potential and corresponding Knapp-Stein integral operator acting between 
induced representations for two consecutive orthogonal Lie algebra on the flat space (or the sphere)
was constructed. Quite remarkably, the residue analysis of analytic continuation of Knapp-Stein 
intertwining map produces residue families on densities, and the factorization identities
correspond to Kummer's relation for Gauss hypergeometric function. We expect that an analogous, though
vector valued version of such identities, will "explain" the origin of spinor residue family operators
and their factorization properties.     
 
Yet another issue is closely related to the spinor version of the holographic deformation 
${\mathcal H}(h)(r)=\sum_{N=1}^\infty{\mathcal M}_{2N}\frac{1}{(N-1)!^2}(\frac{r^2}{4})^{N-1}$ 
of the Yamabe operator $\mathcal{M}_2$, i.e., 
the second order operator on $(M,h)$ obtained as the generating function for the operators 
${\mathcal M}_{2N}$, cf. \cite{Juhl1}. In other words, it is desirable to construct a holographic deformation of 
the Dirac operator. 



It is perhaps worth to note that one of the advantages of the explicit formulas for 
${D}_{2N}^{res}(h; \lambda)$ is that the derivatives of residue families 
with respect to the spectral parameter $\lambda$ yield the $Q$-curvature operators.
It is not clear to the authors, whether there is a reasonable construct for a spinor $Q$-curvature, but 
one can clearly produce the spectral derivatives of spinor residue family operators 
${\slashed D}_{N}^{res}(h; \lambda)$, $N\in{\mathbb N}_0$.


\section*{Appendix A: Spin geometry}\label{Appendix A}

Let $(M^n,h)$ be a semi-Riemannian $Spin$-manifold of signature $(p,q)$, $n=p+q$. 
Then any orthonormal frame $\{s_i\}_i$ fulfills 
$h(s_i,s_j)=\varepsilon_i\delta_{ij}$, where $\varepsilon_i=-1$ 
for $1\leq i\leq p$ and $\varepsilon_i=1$ for $p+1\leq i\leq n$.
 
The Clifford algebra of $(\R^n,\langle\cdot,\cdot\rangle_{p,q})$, 
denoted by $Cl(\R^{p,q})$, is the quotient of tensor algebra of $\R^n$ 
by two sided non-homogeneous ideal generated by relations 
$x\otimes y+y\otimes x=-2\langle x,y\rangle_{p,q}$ for all $x,y\in\R^n$. 

In the even case $n=2m$ the complexified Clifford algebra $Cl_{\mathbb C}(\R^{p,q})$ 
has up to an isomorphism a unique irreducible representation, whereas in the odd case 
$n=2m+1$ it has up to an isomorphism two non-equivalent irreducible representations 
on $\Delta_{p,q}:=\C^{2^m}$. The restriction of this representation to the spin group 
$Spin(p,q)$, regarded as a subgroup of the group of units $Cl^*(\R^{p,q})$, 
is denoted by  $\kappa_{p,q}$. 

The choice of a $Spin$-structure $(Q,f)$ on $(M^n,h)$ gives an associated 
spinor bundle $S(M^n,h):=Q\times_{(Spin_0(p,q),\kappa_{p,q})}\Delta_{p,q}$, 
where $Spin_0(p,q)$ denotes the 
connected component of the spin group containing the identity element. 
Then the Levi-Civita connection 
$\nabla^{h}$ on $(M^n,h)$ lifts to a covariant derivative 
$\nabla^{h,S}$ on the spinor bundle. Furthermore, there is a scalar product $<\cdot,\cdot>$ 
on the spinor bundle, which is parallel with 
respect to $\nabla^{h,S}$ and compatible with the Clifford multiplication:
\begin{align*}
  <X\cdot\psi,\phi>=(-1)^{p+1}<\psi,X\cdot\phi>
\end{align*} 
holds for all $X\in\Gamma(TM)$ and $\psi,\phi\in S(M,h)$. The Dirac operator 
$\slashed{D}=\sum_i\varepsilon_i s_i\cdot \nabla_{s_i}^{h,S}$ acting on $S(M,h)$ is 
formally anti-selfadjoint, 
\begin{align*}
  < \slashed{D}\psi,\phi>_{L^2}=(-1)^p<\psi,\slashed{D}\phi>_{L^2},
\end{align*}
where $<\cdot,\cdot>_{L^2}$ denotes the induced $L^2$-scalar product. 

Let $\widehat{h}:=e^{2\sigma}h$ be a metric conformally related to $h$, $\sigma\in \mathcal{C}^\infty(M)$. 
The spinor bundles for $\widehat{h}$ and $h$ can be identified through a vector bundle isomorphism 
$F_\sigma:S(M,h)\to S(M,\widehat{h})$, and the Dirac operator satisfies the following conformal covariance: 
\begin{align*}
  \widehat{\slashed{D}}(e^{\frac{1-n}{2}\sigma}\widehat{\psi})
     =e^{-\frac{1+n}{2}\sigma}\widehat{\slashed{D}\psi},
\end{align*}
for any $\psi\in\Gamma\big(S(M,h)\big)$ and $\widehat{\cdot}$ denotes evaluation with respect to $\widehat{h}$.  
Conformal odd powers of the Dirac operator 
were constructed in \cite{HS,GMP1, Fischmann}, and are denoted by 
$\mathcal{D}_{2N+1}=\slashed{D}^{2N+1}+LOT$, for $N\in\N_0$ ($N\leq \frac n2$ for even $n$). They satisfy 
\begin{align*}
  \widehat{\mathcal{D}}_{2N+1}(e^{\frac{2N+1-n}{2}\sigma}\widehat{\psi})
     =e^{-\frac{2N+1+n}{2}\sigma}\widehat{\mathcal{D}_{2N+1}\psi},
\end{align*}
for any $\psi\in\Gamma\big(S(M,h)\big)$.


\section*{Appendix B: Poincar\'e-Einstein metric construction}\label{Appendix-Poincare}
Here we briefly review the content of Poincar\'e-Einstein metric construction, \cite{FG3}. Let $(M^n,h)$ be an $n$-dimensional 
semi-Riemannian manifold, $n\geq 3$. On $X:=M\times (0,\varepsilon)$, for $\varepsilon>0$, we 
consider the metric
\begin{align*}
  g_+=r^{-2}(dr^2+h_r),
\end{align*}
for a $1$-parameter family of metrics $h_r$ on $M$ such that $h_0=h$. The requirement of Einstein 
condition on $g_+$ for $n$ odd, 
\begin{align*}
  Ric(g_+)+ng_+=O(r^\infty),
\end{align*}
uniquely determines the family $h_r$, while for $n$ even the conditions 
\begin{align*}
  Ric(g_+)+ng_+=O(r^{n-2}),\\
  \tr(Ric(g_+)+ng_+)=O(r^{n-1}),
\end{align*}
uniquely determine the coefficients $h_{(2)},\ldots,h_{(n-2)}$, $\tilde{h}_{(n))}$ and the trace of $h_{(n)}$ in the 
formal power series
\begin{align*}
  h_r=h+r^2h_{(2)}+\cdots+r^{n-2}h_{(n-2)}+r^n(h_{(n)}+\tilde{h}_{(n)}\log r )+\cdots\quad . 
\end{align*}
For example, we have 
\begin{align*}
  h_{(2)}=-P,\quad h_{(4)}=\frac 14(P^2-\frac{B}{n-4}),
\end{align*}
where $P$ is the Schouten tensor and $B$ is the Bach tensor associated to $h$. 

The metric $g_+$ on $X$ is called Poincar\'e-Einstein metric associated to a semi-Riemannian manifold $(M,h)$. 

All constructions in the present paper, based on the Poincar\'e-Einstein metric, 
depend for even $n$ on the coefficients $h_{(2)},\ldots,h_{(n-2)}$ and $\tr(h_{(n)})$ only. 
Choosing different representatives $h,\widehat{h}\in [h]$ in the conformal class leads to Poincar\'e-Einstein 
metrics $g^1_+$ and $g^2_+$ related by a diffeomorphism $\Phi:U_1\subset X\to U_2\subset X$, 
where both $U_i$, $i=1,2$, contain $M\times\{0\}$, $\Phi|_{M}=\id_M$, and $g_+^1=\Phi^*g^2_+$ (up to a finite 
order in $r$, for even $n$).   


\section*{Appendix C: Gegenbauer polynomials}\label{Appendix}
We summarize several basic conventions and properties 
related to Gegenbauer polynomials used throughout the article.

The Gegenbauer polynomials are defined in terms of their generating function
\begin{align*}
  \frac{1}{(1-2xt+t^{2})^{\alpha}}=\sum_{n=0}^{\infty}C_n^{\alpha}(x) t^{n},
\end{align*}
and satisfy the recurrence relation
\begin{align*}
  C_{n}^\alpha(x) = \frac{1}{n}\left(2x(n+\alpha-1)C_{n-1}^\alpha(x) -
    (n+2\alpha-2)C_{n-2}^\alpha(x)\right)
\end{align*}
with $C_0^\alpha(x)  = 1, C_1^\alpha(x)  = 2 \alpha x$.
Gegenbauer polynomials are solutions of the Gegenbauer differential equation
\begin{align*}
  \left((1-x^{2})\frac{d^2}{dx^2}-(2\alpha+1)x\frac{d}{dx}
     +n(n+2\alpha)\right) y(x)=0,
\end{align*}
and thus can be written as terminating hypergeometric series
\begin{align*}
  C_n^{\alpha}(x)=\frac{(2\alpha)_{{n}}}{n!}
     \,_2F_1\left(\begin{matrix} -n,&2\alpha+n\\ \alpha+\frac{1}{2}\end{matrix};\frac{1-x}{2}\right).
\end{align*}
The explicit form of Gegenbauer polynomials is
\begin{align*}
  C_n^{\alpha}(x)=\sum_{k=0}^{\lfloor n/2\rfloor}
    (-1)^k\frac{\Gamma(n-k+\alpha)}{\Gamma(\alpha)k!(n-2k)!}(2x)^{n-2k},
\end{align*}
whose consequence is the basic formula for derivative of Gegenbauer polynomials 
\begin{align*}
  \frac{d}{dx}C^{\alpha}_{2N}(x)=2\alpha C^{\alpha +1}_{2N-1}(x).
\end{align*}
The even and odd Gegenbauer polynomials are given, in terms of \eqref{eq:DefA} and \eqref{eq:DefB}, by
\begin{align*}
  \frac{(-1)^NN!}{(-\lambda-\frac{n-1}{2})_N}C_{2N}^{-\lambda-\frac{n-1}{2}}(x)
     =&\sum_{j=0}^N (-1)^ja_j^{(N)}(\lambda)x^{2N-2j},\\
  \frac{(-1)^NN!}{2(-\lambda-\frac{n-1}{2})_{N+1}}C_{2N+1}^{-\lambda-\frac{n-1}{2}}(x)
     =&\sum_{j=0}^N (-1)^j b_j^{(N)}(\lambda)x^{2N+1-2j}.
\end{align*}

\vspace{0.5cm}

{\bf Acknowledgment:} It is our pleasure to thank to A. Juhl for his generosity in 
sharing many mathematical ideas with us. The authors gratefully acknowledge the support of the grant GA CR P201/12/G028.

\bibliographystyle{amsalpha}
\bibliography{bibliography}

\newcommand{\etalchar}[1]{$^{#1}$}
\providecommand{\bysame}{\leavevmode\hbox to3em{\hrulefill}\thinspace}
\providecommand{\MR}{\relax\ifhmode\unskip\space\fi MR }
\providecommand{\MRhref}[2]{%
  \href{http://www.ams.org/mathscinet-getitem?mr=#1}{#2}
}
\providecommand{\href}[2]{#2}
\begin{thebibliography}{KKM{\etalchar{+}}78}

\bibitem[BGM05]{BGM}
C.~B{\"a}r, P.~Gauduchon, and A.~Moroianu, \emph{Generalized cylinders in
  semi-{R}iemannian and spin geometry}, Mathematische Zeitschrift \textbf{249}
  (2005), no.~3, 545--580.

\bibitem[Bur93]{Bures}
J.~Bure\v{s}, \emph{Dirac operators on hypersurfaces}, Commentationes
  Mathematicae Universitatis Carolinae \textbf{34} (1993), no.~2, 313--322.

\bibitem[FG11]{FG3}
C.~Fefferman and C.~R. Graham, \emph{The ambient metric ({AM}-178)}, no. 178,
  Princeton University Press, 2011.

\bibitem[Fis13]{Fischmann}
M.~Fischmann, \emph{Conformally covariant differential operators acting on
  spinor bundles and related conformal covariants}, Ph.D. thesis, Humboldt
  Universit\"at zu Berlin, 2013,
  \url{http://edoc.hu-berlin.de/dissertationen/fischmann-matthias-2013-03-04/P%
DF/fischmann.pdf}.

\bibitem[GJMS92]{GJMS}
C.R. Graham, R.W. Jenne, L.~Mason, and G.~Sparling, \emph{Conformally invariant
  powers of the {L}aplacian, {I}: {E}xistence}, Journal of the London
  Mathematical Society \textbf{2} (1992), no.~3, 557.

\bibitem[GMP10]{GMP}
C.~Guillarmou, S.~Moroianu, and J.~Park, \emph{Eta invariant and {S}elberg zeta
  function of odd type over convex co-compact hyperbolic manifolds}, Advances
  in Mathematics \textbf{225} (2010), no.~5, 2464--2516.

\bibitem[GMP12]{GMP1}
\bysame, \emph{Bergman and {C}alder\'on projectors for {D}irac operators},
  Journal of Geometric Analysis (2012), 1--39,
  \url{http://dx.doi.org/10.1007/s12220-012-9338-9}.

\bibitem[GP03]{GoverPeterson}
A.R. Gover and L.J. Peterson, \emph{Conformally invariant powers of the
  {L}aplacian, {$Q-$}curvature, and tractor calculus}, Communications in
  Mathematical Physics \textbf{235} (2003), no.~2, 339--378.

\bibitem[Gra99]{Graham}
C.R. Graham, \emph{Volume and area renormalizations for conformally compact
  {E}instein metrics}, Proc. of 19th Winter School in Geometry and Physics,
  Srni, Czech Rep. (1999).

\bibitem[GZ03]{GZ}
C.R. Graham and M.~Zworski, \emph{Scattering matrix in conformal geometry},
  Inventiones mathematicae \textbf{152} (2003), 89--118.

\bibitem[Hel70]{Helgason}
S.~Helgason, \emph{A duality for symmetric spaces with applications to group
  representations}, Advances in Mathematics \textbf{5} (1970), no.~1, 1--154.

\bibitem[HS01]{HS}
J.~Holland and G.~Sparling, \emph{Conformally invariant powers of the ambient
  {D}irac operator}, ArXiv e-prints (2001), 1--17,
  \url{http://arxiv.org/abs/math/0112033}.

\bibitem[Juh09]{Juhl}
A.~Juhl, \emph{Families of conformally covariant differential operators,
  {$Q-$}curvature and holography}, Birkh\"auser, 2009.

\bibitem[Juh13]{Juhl1}
\bysame, \emph{Explicit formulas for {GJMS}-operators and {$Q-$}curvatures},
  Geometric and Functional Analysis \textbf{23} (2013), 1278--1370,
  \url{http://arxiv.org/abs/1108.0273}.

\bibitem[KKM{\etalchar{+}}78]{KKMOOT}
M.~Kashiwara, A.~Kowata, K.~Minemura, K.~Okamoto, T.~Oshima, and M.~Tanaka,
  \emph{Eigenfunctions of invariant differential operators on a symmetric
  space}, The Annals of Mathematics \textbf{107} (1978), no.~1, 1--39.

\bibitem[Kob13]{kob}
T.~Kobayashi, \emph{F-method for symmetry breaking operators}, Differential
  Geometry and its Applications (2013), 1--18.

\bibitem[Kos75]{Kosmann}
Y.~Kosmann, \emph{Sur les degr\'es conformes des op\'erateurs diff\'erentiels},
  Comptes rendus hebdomadaires des s\'eances de l'Acad\'emie des Sciences,
  S\'erie A \textbf{280} (1975), 229--232.

\bibitem[K{\O}SS13]{KOSS}
T.~{Kobayashi}, B.~{\O}rsted, P.~Somberg, and V.~Sou\v{c}ek, \emph{{Branching
  laws for {V}erma modules and applications in parabolic geometry. I}}, ArXiv
  e-prints (2013), 1--37, \url{http://arxiv.org/pdf/1305.6040v1.pdf}.

\bibitem[PS95]{ps}
P.~Paule and M.~Schorn, \emph{A mathematica version of {Z}eilberger's algorithm
  for proving binomial coefficient identities}, Journal of symbolic computation
  \textbf{20} (1995), no.~5-6, 673--698.

\bibitem[RC01]{cam}
E.~Pedon R.~Camporesi, \emph{Harmonic analysis for spinors on real hyperbolic
  spaces}, Colloquium Mathematicum \textbf{87} (2001), no.~2, 245--286,
  \url{http://calvino.polito.it/~camporesi/Coll2001.pdf}.

\bibitem[Sla66]{sl}
L.~Slater, \emph{{Generalized hypergeometric functions}}, Cambridge University
  Press, 1966.

\end{thebibliography}

\vspace{0.5cm}

Matthias Fischmann, Petr Somberg

Eduard \v{C}ech Institute and Mathematical Institute of Charles University,

Sokolovsk\'a 83, Praha 8 - Karl\'{\i}n, Czech Republic, 

E-mail: fischmann@karlin.mff.cuni.cz, somberg@karlin.mff.cuni.cz.

\end{document}